\documentclass[a4paper,11pt]{amsart}
\usepackage{graphicx}
\usepackage{amsmath,amssymb,amsthm,amsfonts}
\usepackage{amssymb}
\usepackage[active]{srcltx}
\usepackage[colorlinks=true,urlcolor=blue,
citecolor=red,linkcolor=blue,linktocpage,pdfpagelabels,
bookmarksnumbered,bookmarksopen]{hyperref}
\usepackage{amsmath,amssymb,amsthm,amsfonts}
\usepackage{amssymb}

\usepackage[english]{babel}
\usepackage[colorlinks=true]{hyperref}
\usepackage{xcolor}
\usepackage[colorinlistoftodos]{todonotes}

\newtheorem{thm}{Theorem}[section]

\newtheorem{lem}[thm]{Lemma}

\newtheorem{rem}[thm]{Remark}

\newcommand{\R}{{\mathbb{R}}}

\definecolor{darkgreen}{rgb}{0.09, 0.45, 0.27}

\renewcommand{\le}{\leqslant}
\renewcommand{\ge}{\geqslant}

\renewcommand{\l }{\lambda}
\newcommand{\n }{\nabla }

\renewcommand{\L}{\Lambda}

\newcommand{\I}{\mathcal{I}}
\newcommand{\E}{\mathcal{E}}

\newcommand{\Sl}{\Sigma_\lambda}

\newcommand{\beq}{\begin{equation}}
\newcommand{\eeq}{\end{equation}}
\newcommand{\beqs}{\begin{equation*}}
\newcommand{\eeqs}{\end{equation*}}
\newcommand{\beqn}{\begin{eqnarray}}
\newcommand{\eeqn}{\end{eqnarray}}
\newcommand{\beqns}{\begin{eqnarray*}}
\newcommand{\eeqns}{\end{eqnarray*}}
\newcommand{\bdoc}{\begin{document}}
\newcommand{\edoc}{\end{document}}
\newcommand{\be}{\begin{enumerate}}
\newcommand{\ee}{\end{enumerate}}
\newcommand{\bdescr}{\begin{description}}
\newcommand{\edescr}{\end{description}}
\newcommand{\ba}{\begin{array}}
\newcommand{\ea}{\end{array}}
\newcommand{\intR}{\int_{\mathbb R^n}}
\newcommand{\e}{\varepsilon}
 \renewcommand{\(}{\left(}
\renewcommand{\)}{\right)}
\renewcommand{\[}{\left[}
\renewcommand{\]}{\right]}
\newenvironment{Proof}{\noindent{\bf Proof}}{\hfill$\Box$\\[2mm]}

\usepackage{todonotes}

\everymath{\displaystyle}
\numberwithin{equation}{section}

\begin{document}


\parindent 0pc
\parskip 6pt
\overfullrule=0pt

\title[Critical quasilinear elliptic systems]{Radial symmetry of positive solutions to quasilinear Hardy-Sobolev doubly critical systems}

\author {Laura Baldelli$^{\#,+}$, Francesco Esposito$^{*}$, Rafael Lopez-Soriano$^\#$, Berardino Sciunzi$^{*}$}

\date{\today}

\date{\today}

\address{*Dipartimento di Matematica e Informatica, UNICAL,
Ponte Pietro  Bucci 31B, 87036 Arcavacata di Rende, Cosenza, Italy.}

 \email{francesco.esposito@unical.it, berardino.sciunzi@unical.it}
 
 \address{$\#$ Departamento de An\'alisis Matem\'atico, Universidad de Granada, Campus Fuentenueva, 18071 Granada, Spain.}
 
 \email{labaldelli@ugr.es, ralopezs@ugr.es}

 \address{$+$ Institute for Analysis, Karlsruhe Institute of Technology (KIT), D-76128 Karlsruhe, Germany.}
 
 \email{laura.baldelli@kit.edu}


\thanks{L. Baldelli, F. Esposito and R. López-Soriano are partially supported by the PID Project (Spain) PID2024-155314NB-I00. F. Esposito and B. Sciunzi are partially supported by PRIN 2022 project P2022YFAJH {\em Linear and Nonlinear PDE's: New directions and Applications}. L. Baldelli, F. Esposito and B. Sciunzi are members of INdAM. { L. Baldelli and }F. Esposito are partially supported by INdAM - GNAMPA Project 2025 - CUP $E5324001950001$. {Most of the paper was written while L.B. was affiliated with the University of Granada. L. Baldelli is partially supported by the ``Maria de Maeztu'' Excellence Unit IMAG, reference CEX2020-001105-M, funded by MCIN/AEI/10.13039/501100011033/ and by the Deutsche Forschungsgemeinschaft (DFG, German Research Foundation) - Project-ID 258734477 - SFB 1173}}



\begin{abstract}
	The aim of this paper is to prove radial symmetry results for positive weak solutions with finite energy to the following quasilinear doubly critical system 
	\begin{equation*}
		\begin{cases}
			\displaystyle -\Delta_p u\,=\gamma \frac{u^{p-1}}{|x|^p} + u^{p^*-1}+ \nu \alpha u^{\alpha-1} v^\beta & \text{in}\quad\R^n \vspace{0.2cm} \\
			\displaystyle -\Delta_p v\,=\gamma \frac{v^{p-1}}{|x|^p} + v^{p^*-1}+ \nu \beta u^\alpha v^{\beta-1} & \text{in}\quad\R^n,
		\end{cases}
	\end{equation*}
	where $1<p<n$, $\gamma \in [0, \Lambda_{n,p})$ with $\Lambda_{n,p} = \left[(n-p)/p\right]^p$, $\alpha, \beta > 1$ such that $\alpha + \beta = p^*=np/(n-p)$ and $\nu>0$.  
\end{abstract}

\date{\today}

\maketitle

\section{Introduction}

In this paper, we investigate qualitative properties of solutions to the following doubly critical system involving the $p$-Laplace operator:
\begin{equation}\tag{$\mathcal{S}^*$}\label{P}
\begin{cases}
\displaystyle -\Delta_p u\,=\gamma \frac{u^{p-1}}{|x|^p} + u^{p^*-1}+ \nu \alpha u^{\alpha-1} v^\beta & \text{in}\quad\R^n \vspace{0.2cm} \\
\displaystyle -\Delta_p v\,=\gamma \frac{v^{p-1}}{|x|^p} + v^{p^*-1}+ \nu \beta u^\alpha v^{\beta-1} & \text{in}\quad\R^n \vspace{0.2cm} \\
u,v > 0 &  \text{in}\quad\R^n \setminus \{0\}\\
u,v \in \mathcal D^{1,p}(\R^n),
\end{cases}
\end{equation}
where $\gamma \in [0, \Lambda_{n,p})$ with $\Lambda_{n,p} = \left[(n-p)/p\right]^p$ representing the optimal constant in Hardy's inequality for $n > p$. We note that this system is well posed once one assumes that both
$\alpha, \beta > 1$ are real parameters that satisfy
\begin{equation*}
\alpha + \beta = p^*. 
\end{equation*}
Here, $p^* = np/(n-p)$ denotes the critical Sobolev exponent, $\nu > 0$ is the coupling parameter. By the previous homogeneity condition we deduce that $p^*> 2$. Thus, here and in all the paper we will assume that
$$
\frac{2n}{n+2} < p< n.
$$
We point out that $\mathcal D^{1, p}(\mathbb R^n)$ is the completion of $\mathcal C^\infty_c(\mathbb R^n)$, the space of smooth functions with compact support, with respect to the norm 
$$\|u\|:=\left(\intR|\nabla u|^p \, dx \right)^{\frac 1p}.$$
It is well known, by standard regularity theory (see e.g.~\cite{Di,T}), it follows that solutions to \eqref{P} are of class $\mathcal C^{1,\alpha}$ far from the origin.\\

Our main objective is to prove radial symmetry of solutions to system \eqref{P}. When $\gamma = \nu = 0$, the system reduces to the standard $p$-Sobolev critical equation
\begin{equation}\label{eq:Sobolev}
	\begin{cases}
		\displaystyle -\Delta_p u\,=  u^{p^*-1}  & \text{in}\quad\R^n\\
		u > 0 &  \text{in}\quad\R^n.
	\end{cases}
\end{equation}
It is well-known that the Aubin-Talenti bubbles, defined as
\begin{equation}\label{AubinTalenti}
\mathcal{V}_{\lambda,x_0}(x) := \left[\frac{\lambda^{\frac{1}{p-1}}\sqrt[p]{n\left(\frac{n-p}{p-1}\right)^{p-1}}}{\lambda^{\frac{p}{p-1}} + |x - x_0|^{\frac{p}{p-1}}}\right]^{\frac{n-p}{p}}
\end{equation}
is a family of positive radial solutions to \eqref{eq:Sobolev}, with $\lambda > 0$, and $x_0 \in \mathbb{R}^n$. These functions achieve equality in the Sobolev inequality in $\mathbb{R}^n$. The situation in case $p=2$ is fully understood. In fact, in \cite{GNN2}, it has been shown that any positive solution $u$ to \eqref{eq:Sobolev} that satisfies the decay assumption $u(x) = \mathcal{O}(1/|x|^m)$ at infinity for $m > 0$ is radially symmetric and decreasing around some point $x_0 \in \mathbb{R}^n$, i.e.~$u(x)=\mathcal{V}_{\lambda,x_0}(x)$ (with $p=2$). The authors used a refined version of the moving planes technique taking into account the lack of compactness in the whole space $\R^n$, developed by themselves in bounded domains in \cite{GNN}. The problem was completely solved in the seminal paper by Caffarelli, Gidas and Spruck \cite{CGS}, where the authors classified all the solutions to \eqref{eq:Sobolev} with $p=2$ and without an apriori finite energy assumption. Exploiting the Kelvin transformation, they showed that the moving planes procedure can start. Finally, they proved that any solution $u \in H^1_{loc}(\R^n)$ to \eqref{eq:Sobolev} is of the form \eqref{AubinTalenti}, hence solutions are unique up to rescaling and translations.

The situation in the case $p \neq 2$ is much more involved, since the approach with Kelvin transformation is not  applicable for this equation. Under the finite energy assumption, i.e.~$u \in \mathcal{D}^{1,p}(\R^n)$, the classification of all the positive solutions to \eqref{eq:Sobolev}  has been in fact an open and challenging problem recently solved in a series of papers by Damascelli, Merch\'an, Montoro and Sciunzi~\cite{DMMS} ($2n/(n+2)<p<2$), V\'etois~\cite{Vet} ($1<p<2$) and Sciunzi~\cite{Sciu16} ($2<p<n$). In particular, in these articles it was proved that any positive solution to \eqref{eq:Sobolev} that belongs to $\mathcal D^{1,p}(\R^n)$ is of the form \eqref{AubinTalenti}. More recently, Ciraolo, Figalli and Roncoroni~\cite{CFR} used a strategy based on integral estimates to extend the classification of positive $\mathcal D^{1,p}$-solutions to a class of anisotropic $p$-Laplace-type equations in convex cones. A similar approach was used in \cite{CatMonRon} to obtain new classification results for positive, weak solutions to \eqref{eq:Sobolev} which are not a priori in $\mathcal D^{1,p}(\R^n)$. In particular, in \cite{CatMonRon}, the authors managed to obtain the complete classification of positive, weak solutions to \eqref{eq:Sobolev} in the case where $n=2$ and $1<p<2$ or $n=3$ and $3/2<p<2$. The method was recently improved by \cite{Ou} who managed to extend this result to the case where $n\geq 3$ and $p\geq (n+1)/3$ and by Vet\'ois~\cite{Vet2} for $n\geq 4$ and $p>p_n$ for some number $p_n\in (n/3,(n+1)/3)$ such that $p_n\sim n/3+ 1/n$.

Now, we focus our attention on the case $\gamma \neq 0$ and $\nu = 0$; in this setting, system \eqref{P} reduces to the Hardy-Sobolev doubly critical equation
\begin{equation}\label{eq:HardySobolev}
\begin{cases}
-\Delta_p u = \gamma \frac{u^{p-1}}{|x|^p} + u^{p^*-1} & \text{in } \mathbb{R}^n \\
u > 0 & \text{in } \mathbb{R}^n \setminus \{0\}.
\end{cases} 
\end{equation}
Let us now turn to the case $p=2$ so that $\Lambda_{n,2}$ is the best constant in the Hardy-Sobolev inequality for $p=2$. The main contribution in the semilinear case is due to Terracini in \cite{terracini}, where the author, by means of variational arguments and the concentration compactness principle, showed the existence of solutions to a more general version of  the equation \eqref{eq:HardySobolev}. Making use of the Kelvin transformation and of the moving plane technique, she was able to prove that any positive solution with finite energy to \eqref{eq:HardySobolev} is radially symmetric about the origin. Finally, thanks to a detailed ODE's analysis, she gave a complete classification of the solutions to \eqref{eq:HardySobolev}, given by
\begin{equation}\label{eq:terraciniane}
	\mathcal{U}_\lambda(x)=\lambda^{\frac{2-n}{2}}\mathcal{U}\left(\frac{x}{\lambda}\right) \qquad \mbox{ with } \qquad \mathcal{U}(x)= \dfrac{A(n,\gamma)}{|x|^{\mu}\left(1+|x|^{2-\frac{4\mu}{n-2}}\right)^{\frac{n-2}{2}}},
\end{equation}
where 
\begin{equation*}
	\mu:=\frac{n-2}{2}-\sqrt{\left( \frac{n-2}{2}\right)^2-\gamma}, \qquad A(n,\gamma):=\left(\frac{n(n-2-2\mu)^2}{n-2}\right)^{\frac{n-2}{4}},
\end{equation*}
and $\lambda>0$ is a scaling factor. Obviously, $\mathcal{U}_\lambda=\mathcal{V}_{\lambda,0}$ if $\gamma=0$.

The case of $p$-Lapalce operator was firstly treated in \cite{abdellaoui}, where the authors studied the existence of solutions to \eqref{eq:HardySobolev} and a very fine ODE's analysis was carried out. The radial symmetry of the solutions there was an assumption. In \cite{OSV} the authors proved that each solution to \eqref{eq:HardySobolev} with finite energy is radially symmetric (and radially decreasing) about the origin. Having in force the radial symmetry of the solution  it is easy to derive the ordinary differential equation associated to the PDE and fulfilled by the solution $u = u(r)$ and, hence, to apply the results in \cite{abdellaoui}. In particular, they proved also uniqueness up to scaling of the \emph{radial} solutions as well as  the asymptotic behavior are proved showing in particular that, a \emph{radial} solution  \eqref{eq:HardySobolev}, satisfy the following properties:
\begin{equation*}
    \begin{split}
        &\mathcal{U}_{p,\lambda}(x)= \lambda^\frac{ p-n}{p} \mathcal{U}_p\left(\frac{x}{\lambda}\right) \quad \text{ with } \quad \mathcal U_p(x)= u(r)\\
        &\lim_{r\to0^+} r^{\mu_1} \mathcal U_p(r)=C_1,\qquad \lim_{r\to+\infty}r^{\mu_2} \mathcal U_p(r)=C_2,\\
        &\lim_{r\to0^+} r^{\mu_1+1} |\mathcal U_p'(r)|=C_1\mu_1,\qquad \lim_{r\to+\infty}r^{\mu_2+1}|\mathcal{U}_p'(r)|=C_2\mu_2,
    \end{split}
\end{equation*}
for some positive constants $C_1, C_2$. {Obviously, $\mathcal{U}_{p,\lambda}=\mathcal{U}_{\lambda}$ and $\mathcal{U}_{p}=\mathcal{U}$ if $p=2$.} Here and in the sequel $\mu_1, \mu_2\in [0, +\infty)$, $\mu_1<\mu_2$ are two roots of the equation
\begin{equation*}
\mu^{p-2}\left[(p-1)\mu^2-(n-p)\mu\right]+\gamma=0.
\end{equation*}
We remark (for later use) that
\begin{equation*}\label{relgammaj}
0\le \mu_1<\frac{n-p}{p}<\mu_2\le \frac{n-p}{p-1}.
\end{equation*}
Note that when $p=2$ then $\mu_{1}=\tau$ and $\mu_2=n-2-\tau$. Instead, when $p\neq 2$ but $\gamma=0$ then $\mu_1=0$ and $\mu_2=\frac{n-p}{p-1}$.\\

We now turn our attention to the study of \eqref{P}. Nonlinear Schrödinger problems, particularly those of Gross--Pitaevskii type, seems to have strong connections to various physical phenomena. These kind of systems appear naturally in the Hartree--Fock theory for double condensates specifically in binary mixtures of Bose Einstein condensates occupying different hyperfine states while spatially overlapping (see \cite{Esry} and \cite{Frantz} for comprehensive details). 
The solitary-wave solutions to coupled Gross--Pitaevski equations satisfy the following system:
\begin{equation}\label{eq:pBSsystem}
	\begin{cases}
		-\Delta u + V(x) u = \mu u^{2q-1} + \nu u^{q-1} v^{q} & \text{in } \mathbb{R}^n \\
		-\Delta v + V(x) v = \mu v^{2q-1} + \nu u^{q}v^{q-1} & \text{in } \mathbb{R}^n,
	\end{cases}
\end{equation}
where $V$ represents the potential of the system and $1 < q \leq \frac{2^*}{2}$. This formulation is commonly referred to as the Bose--Einstein condensate system. For subcritical regimes, we refer to \cite{ambrosetti}, \cite{bartsch2}, \cite{lin1}, \cite{sirakov}, and \cite{nicola2} for existence and multiplicity results under various assumptions on $V$ and $\nu$. 

In the critical case where $q = \frac{2^*}{2}$, if $V$ is a non-zero constant, system \eqref{eq:pBSsystem} admits only the trivial solution $(0,0)$, that is a consequence of the Pohozaev identity. When $V=0$, \cite{Wang} demonstrated the uniqueness of ground states under specific parameter conditions for generalized systems. Conversely, \cite{pistoia} examined the competitive setting ($\nu < 0$), establishing the existence of infinitely many fully nontrivial solutions that are not conformally equivalent. 

In \cite{ELS}, the authors treated the case of Hardy-type potential $V = -\frac{\gamma}{|x|^2}$. With this choice of $V$, the authors were able to study problem \eqref{P} when $p=2$. What they have proved in \cite{ELS} is that any $(u,v) \in \mathcal{D}^{1,2}(\mathbb{R}^n) \times \mathcal{D}^{1,2}(\mathbb{R}^n)$ solution to \eqref{P} is of synchronized type:
\begin{equation*}
(u,v)=\left(c_1 \, \mathcal{U}_{\lambda_0}, c_2 \, \mathcal{U}_{\lambda_0}\right),
\end{equation*}
where $\mathcal{U}_\lambda$ is given in \eqref{eq:terraciniane}, $\lambda_0>0$ and $c_1,c_2$ are any positive constants satisfying the system
\begin{equation*}
		\left\{\begin{array}{ll}
			\displaystyle{c_1^{2^*-2}+\nu \alpha c_1^{\alpha-2}c_2^\beta}=1  & \vspace{.3cm}\\
			\displaystyle{c_2^{2^*-2}+\nu \beta c_1^{\alpha}c_2^{\beta-2}=1}. &
		\end{array}\right.
\end{equation*}
This result was new also in the case  $\gamma=0$.  Under that assumption the explicit solutions in \eqref{eq:terraciniane} reduce to those of \eqref{AubinTalenti}.\\ 

Doubly critical problems have received significant attention in recent years. The pioneering work of \cite{AFP} investigated general Hardy-Sobolev type systems making use of variational techniques. In cooperative regimes ($\nu > 0$), they established the existence of ground and bound states contingent on parameters $\alpha$, $\beta$, dimension $n$, and a potential function $h$ in the coupling term. We also refer to \cite{ChenZou3, EduRafaAle} for further existence results in the critical regime.

The primary objective of this paper is to show that any positive solutions, with finite energy, to problem \eqref{P} is radially symmetric and radially decreasing around the origin, in the same spirit of \cite{ELS}.

Thus our main result is given by the following:
\begin{thm}\label{thm:symmetry}
	Let $(u,v) \in \mathcal{D}^{1,p}(\mathbb{R}^n) \times \mathcal{D}^{1,p} (\mathbb{R}^n)$ be a solution to \eqref{P}. Then $(u,v)$ is radially symmetric and radially decreasing around the origin.
\end{thm}

\begin{rem}
This result holds (and is new)  also in the case $\gamma=0$.  
\end{rem}

The proof of the radial symmetry of the solutions is based on a fine adaptation of the moving plane procedure for quasilinear elliptic problems. This procedure goes back to the seminal papers of Alexandrov and Serrin \cite{A, BN, GNN, serrin}. The situation in the quasilinear case is much more involved, and in the same spirit of \cite{OSV, Sciu16} we first recover asymptotic estimates for the solutions and its gradient. Part of this estimate are due to \cite{KL20}. The missing part is proved in this paper, and it was inspired by \cite{EMSV}. 
Regarding elliptic systems, the moving plane technique was adapted by Troy in \cite{troy}, where the cooperative case is analyzed. The procedure was also applied for semilinear systems in the half-space in \cite{Dan} and in the whole space by Busca and Sirakov in \cite{busca}. There are several results on the qualitative properties of solutions for cooperative elliptic systems; in particular, we refer to \cite{esposito, ELS} for a complete analysis in the semilinear case and to \cite{BEMV, EMM} for the quasilinear case.

\section{Preliminaries and asymptotic behavior of solutions}

{\bf Notation.} Generic fixed and numerical constants will be denoted by $C$ (with subscript in some case), which may vary line by line, even within a single formula. By $|A|$ we will denote the Lebesgue measure of a measurable set $A$.

\subsection{Preliminary results} 

The aim of this section is to give some preliminary technical lemmas and theorems that will be crucial in the proof of the main result.

The first result, for which we refer to \cite{NS}, summarizes some Caffarelli-Kohn-Nirenberg type inequalities.
\begin{thm}[Proposition 1.1 of \cite{NS}]\label{hardyns}
	Let $r_1\ge 1$, $r_2>0$, $\bar\kappa, \tilde \kappa \in \mathbb R$ such that
	\begin{equation*}\label{rel}
		\frac{1}{r_2}+\frac{\tilde \kappa}{n}=\frac{1}{r_1}+\frac{\bar\kappa-1}{n},
	\end{equation*}
	and with $$0\le \bar\kappa-\tilde \kappa\le 1.$$
    If $\frac{1}{r_2} +\frac{\tilde \kappa}{n}>0$ then for any $u\in \mathcal C^1_c(\mathbb{R}^n)$ it holds
	\begin{equation}\label{ns1}
    \left(\int_{ \mathbb R^n }|x|^{\tilde \kappa r_2}|u|^{r_2} dx\right)^{\frac{1}{r_2}}\le C \left(\int_{\mathbb{R}^n} |x|^{r_1\bar\kappa}|\nabla u|^{r_1} dx \right)^{\frac{1}{r_1}},
    \end{equation}
	where $C$ is a positive constant independent of $u$.
    If $\frac{1}{r_2} +\frac{\tilde \kappa}{n}<0$, then \eqref{ns1} holds for any $u\in \mathcal C^1_c(\mathbb{R}^n\setminus\{0\})$.
\end{thm}
\begin{rem}\label{remexterior}
    As pointed out in \cite{OSV}, the same result works in the case of $u$ defined in exterior domains with the right decay properties.
\end{rem}
We point out that the natural weight associated to our problem is $|\nabla u|^{p-2}$. However, the presence of Hardy potential forces us to give an estimate from below for an expression which involves this weight. In particular, we state here a result whose proof can be found in  \cite{OSV}. 
\begin{lem}{ \cite[Lemma 2.1]{OSV}}\label{Lemlog}
Let $u, v$ be positive and $\mathcal C^1$-functions in a neighbourhood of some point $x_0\in \mathbb{R}^n$. Then it holds
\begin{equation}\label{stimamag2}
\begin{aligned}
&|\nabla u|^{p-2}\nabla u \cdot \nabla\left(u-\frac{v^p}{u^p} u\right)+|\nabla v|^{p-2}\nabla v \cdot \nabla\left(v-\frac{u^p}{v^p} v\right)
\\
& \ge C_p \min\{v^p,u^p\} \left(|\nabla\log u|+|\nabla \log v|\right)^{p-2}|\nabla\log u-\nabla\log v|^2,
\end{aligned}
\end{equation}
near $x_0$ for some constant $C_p$  depending only on $p$.
\end{lem}

\subsection{Decay and asymptotic behavior}

This section is dedicated to decay estimates for positive weak solutions of \eqref{P} and their gradient, which we will be used in the application of the moving plane method. 
The moving plane procedure is strongly related to suitable comparison principles. When the domain is the whole space,  considering problems with a source term involving the Hardy potential, weak comparison principles are naturally related to the use of Hardy-type inequalities that involve the classical radial weights. 

Now, we state the upper and lower estimates of the solutions and the upper one of their gradient, contained in \cite{KL20}, in our case.
\begin{thm} \label{thm:asymptEst}$($Theorem $1.1$, \cite{KL20}$)$
	Let $(u,v) \in \mathcal{D}^{1,p}(\R^n)\times \mathcal{D}^{1,p}(\R^n)$ be a weak solution of \eqref{P}
    with $1<p<n$, $0 \leq \gamma < \Lambda_{n,p}$, $\alpha+\beta=p^*$. Then there exist positive constants $0 < R_0< 1 < R_1$ depending on $n, p, \gamma, \nu, \alpha,\beta$, $u$ and $v$, such that
    \begin{equation}\label{eq:estat0}
    	\frac{c_1^u}{|x|^{\mu_1}}   \leq u(x) \leq \frac{C_1^u}{|x|^{\mu_1}},\quad \frac{c_1^v}{|x|^{\mu_1}}   \leq v(x) \leq \frac{C_1^v}{|x|^{\mu_1}}  \qquad x \in B_{R_0}\setminus\{0\},
    \end{equation}
	and
	\begin{equation}\label{eq:estatInf}
		\frac{c_2^u}{|x|^{\mu_2}}  \leq u(x) \leq \frac{C_2^u}{|x|^{\mu_2}},\quad \frac{c_2^v}{|x|^{\mu_2}}  \leq v(x) \leq \frac{C_2^v}{|x|^{\mu_2}} \qquad x \in (B_{R_1})^c,
	\end{equation}

where $\mu_1$, $\mu_2$ are the unique solutions, see Remark \ref{rem:rootsHardy}, of
\begin{equation}\label{eq:algebraicHardy}
    \mu^{p-2} [(p-1)\mu^2-(n-p)\mu] + \gamma=0,
\end{equation}
$C_1^u, C_2^u$, $C_1^v, C_2^v$ are positive constants depending on $n,p, \gamma, \nu, \alpha,\beta$ and $u$, $c_1^u, c_1^v$ are positive constants depending on $n,p, \gamma, R_0, \mu_1,\nu, \alpha,\beta$ and $u$ or $v$, respectively, $c_2^u, c_2^v$ are positive constants depending on $n,p, \gamma, R_1, \mu_2,\nu, \alpha,\beta$ and $u$ or $v$, respectively.
\end{thm}

Note that, since the constants in front of Hardy's terms coincide, then the further assumptions $(\mathcal H_2)$ and $(\mathcal H_4)$ in \cite{KL20} trivially hold.

\begin{rem}\label{rem:rootsHardy}
Note that, by $0 \leq \gamma < \Lambda_{n,p}$, then \eqref{eq:algebraicHardy} admits two solutions $\mu_1,\mu_2$ such that
$$0 \leq \mu_1 < \frac{n-p}{p} < \mu_2 \leq \frac{n-p}{p-1},$$
since $(n-p)/p$ is the point of the (negative) minimum of the function related to \eqref{eq:algebraicHardy}.
\end{rem}

\begin{thm}\label{stimagradiente1}$($Theorem $1.3$, \cite{KL20}$)$
    Let $(u,v) \in \mathcal{D}^{1,p}(\R^n)\times \mathcal{D}^{1,p}(\R^n)$ be a weak solution of \eqref{P}
    with $1<p<n$, $0 \leq \gamma < \Lambda_{n,p}$ and $\alpha+\beta=p^*$. Then there exist positive constants $\tilde C^u$, $\tilde C^v$ depending on $n,p, \gamma, \nu, \alpha,\beta$ and $u$ or $v$, respectively, such that
    \begin{equation}\label{eq:estat0grada}
    	|\nabla u(x)| \leq \frac{\tilde  C^u}{|x|^{\mu_1+1}}, \quad |\nabla v(x)| \leq \frac{\tilde  C^v}{|x|^{\mu_1+1}}  \qquad x \in B_{R_0},
    \end{equation}
	and
	\begin{equation}\label{eq:estatInfgrada}
		 |\nabla u(x)| \leq \frac{\tilde C^u}{|x|^{\mu_2+1}}, \quad |\nabla v(x)| \leq \frac{\tilde C^v}{|x|^{\mu_2+1}} \qquad x \in (B_{R_1})^c,
	\end{equation}
where $\mu_1, \mu_2$ are roots of \eqref{eq:algebraicHardy} as in Theorem \ref{thm:asymptEst}, and $0 < R_0< 1 < R_1$ are constants depending on $n, p, \gamma,\nu, \alpha,\beta$, $u$ and $v$.
\end{thm}

We end this section by stating and proving the lower estimate of the gradients of solutions of \eqref{P}, taking inspiration from \cite{EMSV}.

\begin{thm}\label{stimagradiente2}
    Let $(u,v) \in \mathcal{D}^{1,p}(\R^n)\times \mathcal{D}^{1,p}(\R^n)$ be a weak solution of \eqref{P}
    with $1<p<n$, $0 \leq \gamma < \Lambda_{n,p}$ and $\alpha+\beta=p^*$. Then there exist positive constants $\tilde c^u$, $\tilde c^v$ depending on $n,p, \gamma, \nu, \alpha,\beta$ and $u$ or $v$, respectively, such that
    \begin{equation}\label{eq:estat0gradb}
    	|\nabla u(x)| \geq \frac{\tilde  c^u}{|x|^{\mu_1+1}}, \quad |\nabla v(x)| \geq \frac{\tilde  c^v}{|x|^{\mu_1+1}}  \qquad x \in B_{R_2},
    \end{equation}
	and
	\begin{equation}\label{eq:estatInfgradb}
		 |\nabla u(x)| \geq \frac{\tilde C^u}{|x|^{\mu_2+1}}, \quad |\nabla v(x)| \geq \frac{\tilde C^v}{|x|^{\mu_2+1}} \qquad x \in (B_{R_3})^c,
	\end{equation}
where $\mu_1, \mu_2$ are roots of \eqref{eq:algebraicHardy} as in Theorem \ref{thm:asymptEst}, and $0 < R_2< 1 < R_3$ are constants depending on $n, p, \gamma,\nu, \alpha,\beta$, $u$ and $v$.
\end{thm}

Before giving the proof, we have to state a preliminary lemma that will be crucial in the sequel but without proving it since it is the Euclidean case of Theorem $4.1$ in \cite{EMSV} where the Finsler setting is treated (more precisely, instead of the Euclidean norm, a more general function $H^{\circ}$ appears).

\begin{lem}\label{tecnicanew}
	Let $w \in \mathcal C^{1,\alpha}_{loc}(\R^n\setminus \{0\})$ be a positive weak solution of the equation
	\begin{equation*}\label{zioa}
		-\Delta_p w - \frac{\gamma}{|x|^p} w^{p-1} = 0 \qquad \text{in } \R^n\setminus \{0\},
	\end{equation*}
	where $0 \leq \gamma < \Lambda_{n,p}$. Assume that there exist positive constants $C_1, c_1$ such that 
	\begin{equation}\label{a_1}
		\frac{c_1}{|x|^{\mu_1}}\le w(x)\le \frac{C_1}{|x|^{\mu_1}},  \qquad \forall x\in \R^n\setminus \{0\},
	\end{equation}
	and suppose that there exists a positive constant $\hat C_1$ such that 
	\begin{equation}\label{assunzionegradiente}
		|\nabla w(x)|\le \frac{\hat C_1}{|x|^{\mu_1+1}},  \qquad \ \  \forall x\in \R^n\setminus \{0\}, 
	\end{equation}
	then 
	\begin{equation*}\label{v_1}
		w(x)=\frac{\overline c_1}{|x|^{\mu_1}}, \text{ with } \overline{c}_1:=\limsup_{|x|\rightarrow 0} |x|^{\mu_1}w(x),
    \end{equation*}
	On the other hand, suppose that there exist positive constants $\tilde C_1$ and $\tilde c_1$ such that 
	\begin{equation*}\label{a_2}
		\frac{\tilde c_1}{|x|^{\mu_2}}\le w(x)\le \frac{\tilde C_1}{|x|^{\mu_2}},\qquad \forall x\in \R^n\setminus \{0\},
	\end{equation*}
	and suppose that there exists a positive constant $\overline C_1$ such that 
	\begin{equation*}\label{assunzionegradiente2}
		|\nabla w(x)|\le \frac{\overline C_1}{|x|^{\mu_2+1}}, \qquad \ \ \forall x\in \R^n\setminus \{0\}, 
	\end{equation*}
	then  
	\begin{equation*}\label{v_2}
		w(x)=\frac{\check{c_1}}{|x|^{\mu_2}}, 
	\text{	with  } \check{c_1}:=\limsup_{|x|\rightarrow +\infty} |x|^{\mu_2}w(x).
    \end{equation*}
\end{lem}

Now we are ready to prove the lower bound estimates of the gradients of solutions of \eqref{P}.

\begin{proof}[Proof of Theorem \ref{stimagradiente2}]
We will prove \eqref{eq:estatInfgradb} since \eqref{eq:estat0gradb} can be proved in a similar way. Suppose by contradiction that there exist sequences of points $x_m$ and $x_m'$ such that  $|x_m|, |x_m'|\rightarrow +\infty$ and
\begin{equation}\label{bau1}
|x_m|^{\mu_2+1}|\nabla u(x_m)|\rightarrow 0 \qquad \text{for }  m\rightarrow +\infty,
\end{equation}
\begin{equation}\label{bau2}
|x_m'|^{\mu_2+1}|\nabla v(x_m')|\rightarrow 0 \qquad \text{for }  m\rightarrow +\infty.
\end{equation}
Let $R_m:=|x_m|$, $R_m':=|x_m'|$ and consider 
$$w_m(x):=R_m^{\mu_2}u(R_mx), \qquad z_m(x):=(R_m')^{\mu_2}v(R_mx).$$ 
For $0<a<A$ fixed and $m$ sufficiently large, from Theorem \ref{thm:asymptEst}, relabeling the constants, we have 
$$\frac{c}{A^{\mu_2}}\le w_n(x)\le \frac{C}{a^{\mu_2}}, \quad \frac{c}{A^{\mu_2}}\le z_n(x)\le \frac{C}{a^{\mu_2}}\qquad \text{in } \overline{B_{A}\setminus B_{a}}.$$ 
Furthermore, recalling the estimate from above of the gradients of the weak solution $(u, v)$ in Theorem \ref{stimagradiente1}, we get
$$|\nabla w_m(x)|\le \frac{\tilde C^{w_m}}{a^{\mu_2+1}}, \quad |\nabla z_m(x)|\le \frac{\tilde C^{z_m}}{a^{\mu_2+1}} \quad\text{in } \overline{B_{A}\setminus B_{a}}.$$
By \cite{L}, $w_m$ and $z_m$ are also uniformly bounded in $\mathcal C^{1,\eta}(K),$ for $0<\eta<1$ and for any compact set $K\subset B_{A}\setminus B_{a}$. Moreover $$w_m(x)\rightarrow w_{a,A}, \quad z_m(x)\rightarrow z_{a,A} \qquad \text{in } B_A\setminus B_a,$$ 
in the norm $\mathcal C^{1,\eta'}$, for $0<\eta' <\eta.$ 
 Moreover, since $(w_m, z_m)$ weakly solves
\begin{equation*}
\begin{cases}
\displaystyle -\Delta_p w_m-\gamma \frac{w_m^{p-1}}{|x|^p} = R_m^{\mu_2(p-1)+p-\mu_2(p^*-1)}(w_m^{p^*-1}+ \nu \alpha w_m^{\alpha-1} z_m^\beta) & \text{in}\quad B_{A}\setminus B_{a} \vspace{0.2cm} \\
\displaystyle -\Delta_p z_m-\gamma \frac{z_m^{p-1}}{|x|^p} = R_m^{\mu_2(p-1)+p-\mu_2(p^*-1)}(z_m^{p^*-1}+ \nu \beta w_m^\alpha z_m^{\beta-1}) & \text{in}\quad B_{A}\setminus B_{a}\vspace{0.2cm}
\end{cases}
\end{equation*}
and $w_m\rightharpoonup w_{a,A}$, $z_m\rightharpoonup z_{a,A}$ in $\mathcal D^{1,p}(B_{A}\setminus B_{a})$, we deduce that 
\begin{equation*}
\begin{cases}
\displaystyle -\Delta_p w_{a,A}-\gamma \frac{w_{a,A}^{p-1}}{|x|^p} =0 & \text{in}\quad B_{A}\setminus B_{a} \vspace{0.2cm} \\
\displaystyle -\Delta_p z_{a,A}-\gamma \frac{z_{a,A}^{p-1}}{|x|^p} = 0 & \text{in}\quad B_{A}\setminus B_{a}\vspace{0.2cm}
\end{cases}\end{equation*}
Now we take $a_j=1/j$ and $A_j=j$, for $j\in \mathbb N$ and we construct $w_{a_j,A_j}, z_{a_j,A_j}$ as above. For $j$ goes to infinity, using a standard diagonal process, we construct a limiting profile $(w_{\infty},z_\infty)$ so that 
\begin{equation*}
\begin{cases}
\displaystyle -\Delta_p w_{\infty}-\gamma \frac{w_{\infty}^{p-1}}{|x|^p} =0 & \text{in}\quad \R^n\setminus \{0\} \vspace{0.2cm} \\
\displaystyle -\Delta_p z_{\infty}-\gamma \frac{z_{\infty}^{p-1}}{|x|^p} = 0 & \text{in}\quad \R^n\setminus \{0\}\vspace{0.2cm}
\end{cases}\end{equation*}
with $w_{\infty}\equiv w_{a_j,A_j}$ and $z_{\infty}\equiv z_{a_j,A_j}$ in $B_{A_j}\setminus B_{a_j}$.

Since both $w_{\infty}$ and $z_{\infty}$ satisfy the assumptions \eqref{a_1} and \eqref{assunzionegradiente}
of the Lemma \ref{tecnicanew}, we get 
\begin{equation*}
    w_{\infty}(x)= \frac{\overline c_1}{|x|^{\mu_2}},\qquad  z_{\infty}(x)= \frac{\overline c_2}{|x|^{\mu_2}},
\end{equation*}
where $\overline{c}_1:=\limsup_{|x|\rightarrow 0} |x|^{\mu_2}w_{\infty}(x)$, $\overline{c}_2:=\limsup_{|x|\rightarrow 0} |x|^{\mu_2}z_{\infty}(x)$.

Now we set $y_m=x_m/R_m$ and $y_m'=x_m'/R_m$, and by \eqref{bau1} and \eqref{bau2}, we deduce that $|\nabla w_m(y_m)|$ and $|\nabla z_m(y_m')|$ tend to zero as $R_m$ tends to infinity. This fact and the uniform convergence of the gradients imply that there exist $\overline y, \overline y'\in \partial B_{1}$ such that $y_m\to\overline y$, $y_m'\to\overline y'$
$$|\nabla w_{\infty}(\overline y)|=0=|\nabla z_{\infty}(\overline y')|.$$
This is a contradiction since the functions $w_{\infty}$ and $z_{\infty}$ have no critical points.
\end{proof}

\section{Radial symmetry of the solutions}

To prove the symmetry and monotonicity result, we need a fine version of the moving plane technique that is inspired by \cite{OSV}. To exploit it we need the following notations. We will study the symmetry of the solutions in the $\eta$-direction for any $\eta\in \mathbb S^{{n-1}}$ (i.e. $|\eta|=1$).
Since the problem is invariant up to rotations we fix $\eta=e_1$ and we let
$$T_\lambda=\left\{x\in \R\,:\, x_1=\lambda\right\},$$ $$\Sigma_\lambda=\left\{x\in \R\,:\, x_1<\lambda\right\},$$ $$x_\lambda= R_\lambda(x)=(2\lambda-x_1, x')\in \mathbb R\times \mathbb R^{n-1},$$ $$u_\lambda(x)=u(x_\lambda),$$ $$v_\lambda(x)=v(x_\lambda).$$ 

We are now able to prove our symmetry result.

First of all we underline that it is easy to see that $(u_\l,v_\lambda)$ solves
\begin{equation}\tag{$\mathcal{S}^*_\lambda$}\label{pbriflessa}
\begin{cases}
-\Delta_p u_\lambda=\frac{\gamma}{|x_\lambda|^p}u_\lambda^{p-1}+u_\lambda^{p^*-1} + \nu\alpha u_\lambda^{\alpha-1} v_\lambda^{\beta} & \text{in}\quad\R^n \vspace{0.2cm} \\
-\Delta_p v_\lambda=\frac{\gamma}{|x_\lambda|^p}v_\lambda^{p-1}+v_\lambda^{p^*-1} + \nu\beta u_\lambda^{\alpha} v_\lambda^{\beta-1} & \text{in}\quad\R^n \vspace{0.2cm}\\
u_\lambda, v_\lambda > 0 \quad & \text{in}\quad\R^n .
\end{cases}
\end{equation} 

In what follows we set
$$\Lambda_0:=\left\{\lambda<0\,:\, u\le u_t, v \le v_t\,\, \mbox{in}\,\, \Sigma_t, \,\, \forall\, t\le \lambda\right\}.$$
If $\Lambda_0\neq \emptyset$ we denote by $\lambda_0:=\sup\Lambda_0$.\\

\begin{proof}[Proof of Theorem \ref{thm:symmetry}]

We prove the result  studying, sometimes in different approaches, the singular case $1<p<2$ and the degenerate case $p>2$ as in \cite{OSV}. For $p=2$ we refer to \cite{ELS}. We divide the proof in two steps. \\ 

{\bf Step 1:} 	$\Lambda_0\neq \emptyset$.

\vspace{0.3cm}

To prove that $\Lambda_0\neq \emptyset$, it is necessary to prove the existence of $\lambda<0$ with $|\lambda|$ sufficiently large such that $u\le u_t$ and $v\le v_t$ in $\Sigma_t$ for every $t\le \lambda$.\\
We denote by $R_0$ and $R_1$ the radii given by \eqref{eq:estat0} and \eqref{eq:estat0grada}, \eqref{eq:estatInf} and \eqref{eq:estatInfgrada}, while by $R_2$ and $R_3$ the radii given by \eqref{eq:estat0gradb}, \eqref{eq:estatInfgradb} for the entire proof,  and we firstly observe that for $|\overline{\l}|>\max(R_1,R_3)$ one has, by  \eqref{eq:estat0} and \eqref{eq:estatInf}, that there exists $\tilde{R}_0:=\tilde{R}_0(\overline{\l})$ such that $\tilde{R}_0<R_0$, $B_{\tilde{R}_0}(0_{\overline{\l}})\subset \Sigma_{\overline{\l}}$ and
\begin{equation*}\label{stimartilde}
\displaystyle \sup_{x\in B_{\tilde{R}_0}(0_{\overline{\l}})} u(x) < \displaystyle \inf_{x\in B_{\tilde{R}_0}(0_{\overline{\l}})}u_{\overline{\l}}(x) \ \ \text{ and } \sup_{x\in B_{\tilde{R}_0}(0_{\overline{\l}})} v(x) < \displaystyle \inf_{x\in B_{\tilde{R}_0}(0_{\overline{\l}})}v_{\overline{\l}}(x).
\end{equation*}
Therefore,  exploiting again \eqref{eq:estat0} and \eqref{eq:estatInf}, we deduce that for every $\lambda < \overline{\l}$ one has
$$\displaystyle \sup_{x\in B_{\tilde{R}_0}(0_\l)} u(x) \leq \displaystyle \inf_{x\in B_{\tilde{R}_0}(0_\l)}u_{\l}(x) \ \ \text{ and } \sup_{x\in B_{\tilde{R}_0}(0_\l)} v(x) \leq \displaystyle \inf_{x\in B_{\tilde{R}_0}(0_\l)}v_{\l}(x),$$
which gives that $u<u_\l$ and $v<v_\l$ in $B_{\tilde{R_0}}(0_\l)\subset\Sl$ for every $\l\le\overline{\l}$ and with $\tilde{R_0}$ independent of $\l$.

Now, let us consider $R>0$ and a cut-off function $\varphi_R\in \mathcal C^\infty_0(B_{2R}(0))$ such that $0\le \varphi_R \le 1$, $\varphi_R \equiv 1$ on $B_R(0)$ and $|\nabla \varphi_R| \le 2/R$. We denote by $\Sigma_\l'=\Sigma_\l\setminus B_{\tilde{R}_0}(0_\l)$ and by $\hat B_{\rho}:=B_\rho(0)\cap \Sigma_\l'$ for $\rho>0$. Finally, we set
\begin{equation*}\label{eq:diff_funct}
	\begin{split}
	\xi_{\lambda}(x):=u(x)-u_\lambda(x) \ \  &\text{ and }\  \ \ \zeta_{\lambda}:=v(x)-v_\lambda(x),\\
	\xi_{p, \lambda}(x):=u^p(x)-u^p_\lambda(x) \ \  &\text{ and }\  \ \ \zeta_{p, \lambda}(x):=v^p(x)-v^p_\lambda(x),
\end{split}
\end{equation*}
and we recall that for every $\lambda<-  \tilde R_0$ we have
$$u,v \in L^{p^*}(\Sigma_\lambda) \cap L^\infty(\Sigma_\lambda) \cap \mathcal C^0(\overline{\Sigma_\lambda}).$$
since $u,v \in \mathcal D^{1,p}(\R^n)$ and by Theorems \ref{thm:asymptEst} and \ref{stimagradiente1}.
 If $\theta>\max\{2, p\}$ and $\l \le \overline{\l}$, we consider
\begin{equation*}\label{test}
\begin{split}
\Phi_{1, \l}=\varphi_R^\theta u^{1-p}\xi_{p, \lambda}^+\chi_{\Sigma_\l},\qquad \Phi_{2, \l}=\varphi_R^\theta u_\l^{1-p}\xi_{p, \lambda}^+\chi_{\Sigma_\l},\\
\Psi_{1, \l}=\varphi_R^\theta v^{1-p}\zeta_{p, \lambda}^+\chi_{\Sigma_\l},\qquad \Psi_{2, \l}=\varphi_R^\theta v_\l^{1-p}\zeta_{p, \lambda}^+\chi_{\Sigma_\l},
\end{split}
\end{equation*}
We remark that ${\rm supp}(\Phi_{j, \l}),  
 \ {\rm supp}(\Psi_{j, \l}) \subset \hat B_{2R}$ for $j=1, 2$.
Hence, we take $\Phi_{1, \l}$ and $\Psi_{1, \l}$ as a test function in both equation of \eqref{P} respectively, while $\Phi_ {2, \l}$ and $\Psi_ {2, \l}$ in \eqref{pbriflessa} respectively. Then, we   subtract the first equation of \eqref{P} with the first of \eqref{pbriflessa}, and we do the same for both the second equations. Hence, one gets
\begin{equation}\label{sim1}
\begin{aligned}
\int_{\hat B_{2R}} &\left(|\nabla u|^{p-2}\nabla u\cdot\nabla\Phi_{1, \l}-|\nabla u_\l|^{p-2}\nabla u_\l\cdot\nabla\Phi_{2, \l}\right)\, dx\\
=&\gamma\int_{\hat B_{2R}}\left(\frac{1}{|x|^p} - \frac{1}{|x_\l|^p}\right)\varphi_R^\theta\xi_{p, \lambda}^+ \, dx + \int_{\hat B_{2R}}(u^{p^*-p}-u_\l^{p^*-p})\varphi_R^\theta\xi_{p, \lambda}^+\, dx\\
&+\nu \alpha \int_{\hat B_{2R}}(u^{\alpha-p}v^\beta-u_\l^{\alpha-p} v_\l^\beta)\varphi_R^\theta\zeta_{p, \lambda}^+ \, dx,
\end{aligned}
\end{equation}
\begin{equation}\label{sim2}
\begin{aligned}
\int_{\hat B_{2R}}&\left(|\nabla v|^{p-2}\nabla v\cdot\nabla\Psi_{1, \l}-|\nabla v_\l|^{p-2}\nabla v_\l\cdot\nabla\Psi_{2, \l}\right) \, dx\\
=&\gamma\int_{\hat B_{2R}}\left(\frac{1}{|x|^p} - \frac{1}{|x_\l|^p}\right)\varphi_R^\theta \zeta_{p, \lambda}^+\,dx +\int_{\hat B_{2R}}(v^{p^*-p}-v_\l^{p^*-p})\varphi_R^\theta\zeta_{p, \lambda}^+ \, dx\\
&+\nu \beta \int_{\hat B_{2R}}(u^\alpha v^{\beta-p}-u_\l^\alpha v_\l^{\beta-p})\varphi_R^\theta\zeta_{p, \lambda}^+ \, dx,
\end{aligned}
\end{equation}
and, since $|x|\ge |x_\l|$ in $\Sl$, one has that the first terms on the right hand side of \eqref{sim1} and \eqref{sim2} are negative. Using the very definition of the test functions, we get
\begin{equation}\label{inizio}
\begin{aligned}
\I_1:=&\int_{\hat B_{2R}}\varphi_R^\theta \left(|\nabla u|^{p-2}\nabla u\cdot\nabla(u^{1-p}\xi_{p,\lambda}^+)-|\nabla u_\l|^{p-2}\nabla u_\l\cdot\nabla(u_\l^{1-p}\xi_{p,\lambda}^+) \right) \, dx\\
\le&-\theta\int_{\hat B_{2R}} \varphi_R^{\theta-1} u^{1-p} \xi_{p,\lambda}^+ |\nabla u|^{p-2}\nabla u\cdot \nabla \varphi_R \, dx \\&+\theta\int_{\hat B_{2R}} \varphi_R^{\theta-1} u_\l^{1-p} \xi_{p,\lambda}^+|\nabla u_\l|^{p-2}\nabla u_\l \cdot \nabla \varphi_R \,dx \\
&+\int_{\hat B_{2R}} \left(u^{p^*-p}-u_\l^{p^*-p} \, dx\right) \varphi_R^\theta\xi_{p,\lambda}^+ \, dx + \nu \alpha \int_{\hat B_{2R}}\left(u^{\alpha-p} v^\beta-u_\l^{\alpha-p} v_\lambda^\beta \right)\varphi_R^\theta\xi_{p,\lambda}^+ \, dx \\
=:& \ \I_2+\I_3+\I_4+\I_5.
\end{aligned}
\end{equation}
Analogously, we obtain a similar estimate on \eqref{sim2}
\begin{equation}\label{iniziosim2}
\begin{aligned}
\E_1:=&\int_{\hat B_{2R}}\varphi_R^\theta \left(|\nabla v|^{p-2}\nabla v\cdot\nabla(v^{1-p}\zeta_{p,\lambda}^+)-|\nabla v_\l|^{p-2}\nabla v_\l\cdot\nabla(v_\l^{1-p}\zeta_{p,\lambda}^+)\right) \, dx\\
\le &- \theta\int_{\hat B_{2R}}\varphi_R^{\theta-1} v^{1-p} \zeta_{p,\lambda}^+ |\nabla v|^{p-2} \nabla v \cdot \nabla \varphi_R \, dx \\&+ \theta \int_{\hat B_{2R}} \varphi_R^{\theta-1}v_\l^{1-p}\zeta_{p,\lambda}^+|\nabla v_\l|^{p-2}\nabla v_\l \cdot \nabla \varphi_R \,dx \\
&+\int_{\hat B_{2R}}\left(v^{p^*-p}-v_\l^{p^*-p}\right) \varphi_R^\theta\zeta_{p,\lambda}^+ \,dx + \nu \beta \int_{\hat B_{2R}}\left(u^\alpha v^{\beta-p}-u_\l^\alpha v_\lambda^{\beta-p}\right)\varphi_R^\theta\zeta_{p,\lambda}^+\,dx \\
=:& \ \E_2+\E_3+\E_4+\E_5.
\end{aligned}
\end{equation}
In the sequel, we will focus our attention only on \eqref{inizio}, since the estimates on \eqref{iniziosim2} will be very similar; we start by estimating $\I_1$. By using \eqref{stimamag2} it yields that for $p>2$ one has
\begin{equation}\label{I1}
\begin{aligned}
\I_1&=\int_{\hat B_{2R}}\varphi_R^\theta \left(|\nabla u|^{p-2}\nabla u\cdot\nabla(u^{1-p}\xi_{p,\lambda}^+)-|\nabla u_\l|^{p-2}\nabla u_\l\cdot \nabla(u_\l^{1-p}\xi_{p,\lambda}^+)\right) \,dx\\
&=\int_{\hat B_{2R}\cap \{u\ge u_\l\}} \varphi_R^\theta \left(|\nabla u|^{p-2}\nabla u\cdot\nabla\left(u - \frac{u_\l^{p}}{u^p}u\right)+|\nabla u_\l|^{p-2}\nabla u_\l\cdot \nabla \left(u_\l - \frac{u^{p}}{u_\l^p}u_\l\right)\right) \, dx\\
&\ge C_p \int_{\hat B_{2R}\cap \{u\ge u_\l\}} \varphi_R^\theta u_\l^p\(|\n \log u|+|\n\log u_\l|\)^{p-2}|\n\log u-\n\log u_\l|^2 \, dx
\\
&\ge C_p\int_{\hat B_{2R}\cap \{u\ge u_\l\}} \varphi_R^\theta\(\frac{u_\l}{u}\)^p u^2\(|\n u|+|\n u_\l|\)^{p-2}|\n\log u -\n\log u_\l|^2 \, dx
\\
&\ge c_1 \int_{\hat B_{2R}\cap \{u\ge u_\l\}} \varphi_R^\theta u^2\(|\n u|+|\n u_\l|\)^{p-2}|\n\log u -\n\log u_\l|^2 \, dx,
\end{aligned}
\end{equation}
while for $1<p<2$ we obtain
\begin{equation}\label{I1bis}
\begin{aligned}
\I_1&\ge C_p \int_{\hat B_{2R}\cap \{u\ge u_\l\}} \varphi_R^\theta  u_\l^p\frac{|\n\log u-\n\log u_\l|^2}{(|\n \log u|+|\n \log u_\l|)^{2-p}} \,dx
\\
&\ge C_p\int_{\hat B_{2R}\cap \{u\ge u_\l\}} \varphi_R^\theta u_\l^2 \frac{|\n\log u -\n\log u_\l|^2}{\(|\n u|+|\n u_\l|\)^{2-p}} \, dx.
\end{aligned}
\end{equation}
We point out that in the last line of \eqref{I1} we used that
\begin{equation}\label{stimaulu}
\frac{u_\l}{u}\ge \tilde c \quad \mbox{in}\ \Sigma_\l,
\end{equation}
and $c_1:=C_p\tilde{c}^p$.
Indeed if $x \in \Sl\setminus B_{R_1}(0_\lambda)$ then from \eqref{eq:estat0} and \eqref{eq:estatInf} one has (recall that $|x|\ge |x_\l|$)
\begin{equation*}
	\frac{u_\l}{u}\ge \tilde c_1\frac{|x|^{\mu_2}}{|x_\l|^{\mu_2}}\ge \tilde c_1.
\end{equation*}
where $\tilde c_1=c_2^{u_\lambda}/C_2^{u_\lambda}$. Otherwise, if $x \in  \Sl\cap B_{R_1}(0_\lambda)$ then
\begin{equation*}
\frac{u_\l}{u}\ge \tilde c_1|\overline{\l}|^{\mu_2} \inf_{x\in B_{R_1}(0)} u(x) \ge \tilde c_2,
\end{equation*}
where we take $\tilde c=\min(\tilde c_1,\tilde c_2)$ and $|\bar \lambda|>R_1$.
Now, using Schwarz inequality, \eqref{eq:estatInf} and \eqref{eq:estatInfgrada}, we deduce that
\begin{equation}\label{I2finale}
\begin{aligned} \I_2 &=-\theta\int_{\hat B_{2R}} \varphi_R^{\theta-1} u^{1-p} \xi_{p,\lambda}^+ |\nabla u|^{p-2}\nabla u\cdot \nabla \varphi_R \, dx \\
&\le \theta \int_{\hat B_{2R}\cap \{u\ge u_\l\}} \varphi_R^{\theta-1} u\left(1-\(\frac{u_\l}{u}\)^p\right)|\nabla u|^{p-1}|\n\varphi_R| \,dx\\
&\le \frac{2\theta}{R}\int_{\hat B_{2R}\setminus \hat B_R}u|\n u|^{p-1} \,dx  \le \frac{C}{R}\int_{\hat B_{2R}\setminus \hat B_{R}}\frac{1}{|x|^{(\mu_2+1)(p-1)+\mu_2}} \,dx \le \frac{C}{R^\kappa}
\end{aligned}
\end{equation}
where  
\begin{equation}\label{kappa}
\kappa:=p\mu_2+p-n
\end{equation}
which is strictly positive since $\mu_2>\frac{n-p}{p}$. Now we want to estimate $\I_3$; using \eqref{eq:estatInf} and \eqref{stimaulu}, we are able to deduce that
\begin{equation}\label{I3finale}
\begin{aligned} \I_3&= \theta \int_{\hat B_{2R}} \varphi_R^{\theta-1} u_\l^{1-p} \xi_{p,\lambda}^+|\nabla u_\l|^{p-2}\nabla u_\l \cdot \nabla \varphi_R \, dx \\
&\le C\int_{\hat B_{2R}\setminus \hat B_R} \theta \varphi_R^{\theta-1}u_\l^{1-p}\(u^p-u_\l^p\)^+|\n u_\l|^{p-1}|\n \varphi_R| \, dx
\displaystyle
\\
&\le \frac 2 R \int_{(\hat B_{2R}\setminus \hat B_R) \cap \{u\ge u_\l\}}u_\l\(\(\frac{u}{u_\l}\)^p-1\)|\n u_\l|^{p-1} \, dx
\\
&\le \frac 2 R \int_{(\hat B_{2R}\setminus \hat B_R) \cap \{u\ge u_\l\}}u_\l \(\frac{u}{u_\l}\)^p |\n u_\l|^{p-1}\, dx \le \frac C R  \int_{(\hat B_{2R}\setminus \hat B_R) \cap \{u\ge u_\l\}}u|\n u_\l|^{p-1} \, dx\\
&\le \frac C R  \(\int_{ \hat B_{2R}\setminus \hat B_R} |\n u_\l|^p \,dx \)^{\frac{p-1}{p}}\(\int_{\hat B_{2R}\setminus \hat B_R}u^p \,dx\)^{\frac 1 p} \,dx \\&\le \frac C R \(\int_{\hat B_{2R}\setminus \hat B_R}\frac{1}{|x|^{\mu_2p}} \, dx\)^{\frac 1 p}\le \frac {C}{ R^\frac{\kappa}{p}}.
\end{aligned}
\end{equation}
For the term $\I_4$ we first note that  
\begin{equation*}
    \begin{aligned}
\I_4&= \int_{\hat B_{2R} \cap \{u\ge u_\l\} }\left(u^{p^*-p}-u_\l^{p^*-p}\right) \varphi_R^\theta \xi_{p,\lambda}^+ \, dx = \int_{\hat B_{2R} \cap \{u\ge u_\l\} }\left(\frac{u^{p^*-1}}{u^{p-1}}-\frac{u_\l^{p^*-1}}{u_\l^{p-1}}\right) \varphi_R^\theta \xi_{p,\lambda}^+ \, dx\\
&\le  \int_{\hat B_{2R} \cap \{u\ge u_\l\} }\frac{1}{u^{p-1}}\left(u^{p^*-1}-u_\l^{p^*-1}\right)\varphi_R^\theta \xi_{p,\lambda}^+ \, dx,
    \end{aligned}
\end{equation*}
recalling that $p^* \geq 2$, then applying twice Mean Value Theorem and using \eqref{eq:estatInf} we deduce that
\begin{equation}\label{I4}
\I_4 \le c_p\int_{\hat B_{2R}}  u^{p^*-2} \varphi_R^\theta(\xi_\lambda^+)^2 \,dx \le c_p \int_{\hat B_{2R}} \frac{1}{|x|^{\mu_2(p^*-2)}}\varphi_R^\theta(\xi_\lambda^+)^2 \,dx.
\end{equation}
Let us consider the real function $g(t)=\log(a+t(b-a))$ where $0<a\leq b$, then making use of Taylor's formula, we get
$$\log b=\log a +(b-a)\int_0^1 \frac{1}{a+t(b-a)},$$ and since $t\in [0, 1]$, we easily obtain
\begin{equation}\label{usare} b-a = \frac{\log b-\log a}{\displaystyle \int_0^1 \frac{1}{a+t(b-a)}}\le b (\log b-\log a).
\end{equation}
Choosing $b=u$ and $a=u_\l$ in \eqref{usare}, and thanks to \eqref{I4} we deduce
\begin{equation*}\label{i42}
\begin{aligned}
	\I_4 &\le c_p \int_{\hat B_{2R} \cap \{u\ge u_\l\}} \frac{1}{|x|^{\mu_2(p^*-2)}}\varphi_R^\theta(u-u_\l)^2 \, dx\\
    &\le C \int_{\hat B_{2R}\cap \{u\ge u_\l\}}\frac{1}{|x|^{\mu_2(p^*-2)}} \varphi_R^\theta u^2\(\log u-\log u_\l\)^2 \, dx
	\\
	&\le C \int_{\hat B_{2R} \cap \{u\ge u_\l\}} \frac{1}{|x|^{\mu_2 p^*}} \varphi_R^\theta \left((\log u-\log u_\l)^+\right)^2 \, dx.\\
\end{aligned}
\end{equation*}
Furthermore, we can rewrite the previous estimate as follows
\begin{equation}\label{I4bis}
\begin{aligned}
	\I_4 &\le C \int_{\hat B_{2R} \cap \{u\ge u_\l\}} \frac{1}{|x|^{\kappa^*-2{ \vartheta}+2}}\left(\varphi_R^{\frac \theta 2}(\log u-\log u_\l)^+\right)^2 \, dx
	\\
	&\le  \frac{C}{|\l|^{\kappa^*}} \int_{\hat B_{2R} \cap \{u\ge u_\l\}}|x|^{2{ \vartheta}-2}\left(\varphi_R^{\frac \theta 2}(\log u-\log u_\l)^+\right)^2 \, dx,
\end{aligned}
\end{equation}
where
$$\kappa^*:=\mu_2(p^*-p)-p,\qquad 2{\vartheta}:=-[(\mu_2+1)(p-2)+2\mu_2].$$
We note that $\kappa^*-2{ \vartheta}+2=\mu_2p^*$, and $\kappa^*>0$ since $\mu_2>\frac{n-p}{p}$. Applying Theorem \ref{hardyns} to the right hand side of \eqref{I4bis}, where $r_1=r_2=2$ which implies that
$$\tilde \kappa:={\ \bar\kappa}-1=-\frac{(\mu_2+1)p}{2}$$
and that
$$\frac 1 2 +\frac{\tilde \kappa}{n} =\frac{n-\mu_2p-p}{2n}<0$$ since $\mu_2>\frac{n-p}{p}$.  Hence, we obtain
\begin{equation}\label{stimaI4}
\I_4 \le \frac{C}{|\l|^{\kappa^*}} \int_{\hat{B}_{2R} \cap \{u\ge u_\l\}}|x|^{2{\vartheta}}|\nabla (\varphi_R^{\frac \theta 2}(\log u -\log u_\l)^+)|^2 \, dx,
\end{equation}
and, in order to give an estimate of the right hand side of \eqref{stimaI4}, we have to distinguish two cases: $p>2$ and $1<p<2$.
Let us consider the first case, $p>2$ in \eqref{stimaI4}. { By using \eqref{eq:estatInfgradb} and \eqref{stimaulu}, we have}
\begin{equation}\label{stimaI4pmag2}\begin{aligned}
\I_4 \le & \frac{C}{|\l|^{\kappa^*}} \int_{\hat{B}_{2R}\cap\{u\ge u_\l\}} \frac{1}{|x|^{(\mu_2+1)(p-2)}} \varphi_R^\theta u^2|\nabla \log u - \n \log u_\l|^2 \, dx
\\
&+\frac{C}{|\l|^{\kappa^*}} \int_{\hat{B}_{2R}\cap\{u\ge u_\l\}}|x|^{2{\vartheta}}\(\log u - \log u_\l\)^2|\n \varphi_R|^2 \, dx
\\
\le & \frac{C}{|\l|^{\kappa^*}} \int_{\hat{B}_{2R}\cap\{u\ge u_\l\}}\varphi_R^\theta u^2|\n u|^{p-2} |\nabla \log u - \n \log u_\l|^2 \, dx + \frac{C}{|\l|^{\kappa^*}R^2} \int_{\hat{B}_{2R}\setminus\hat{B}_{R}}|x|^{2\theta} \,dx
\\
\le & \frac{C}{|\l|^{\kappa^*}} \int_{\hat{B}_{2R}\cap\{u\ge u_\l\}}\varphi_R^\theta u^2\(|\n u|+|\n u_\l|\)^{p-2} |\nabla \log u - \nabla \log u_\l|^2 \, dx + \frac{C}{|\l|^{\kappa^*}R^\kappa}
\end{aligned}
\end{equation}
where $\kappa$ is defined in \eqref{kappa}.
The evaluation of $\I_5$ is a delicate issue within this argument. Note that, in particular, we could have either $\alpha-p<1$, or $\beta-p<1$, or both. In all these cases we have to take into account the estimate \eqref{stimaulu} which gives us
\begin{equation}\label{eq:ratiobdd}
\tilde{c} \le \frac{u_\lambda}{u} \le 1 \qquad \text{in } \Sigma_{\lambda} \cap \text{supp}(\xi_\lambda^+).
\end{equation}
Now we are ready to estimate the coupling term $\I_5$.
\begin{equation}\label{I5_0}
\begin{aligned} 
\I_5= & \nu \alpha \int_{\hat B_{2R}}\left(u^{\alpha-p}v^\beta-u_\l^{\alpha-p}v_\lambda^\beta\right)\varphi_R^\theta\xi_{p,\lambda}^+ \,dx \\
=& \nu \alpha  \int_{\hat B_{2R} \cap \{u\ge u_\l\}} \left(u^{\alpha-p}v^\beta-u_\l^{\alpha-p} v_\lambda^\beta\right) \varphi_R^\theta (u^p-u_\l^p) \,dx\\
\le & \nu \alpha \int_{\hat B_{2R} \cap \{u\ge u_\l\}} \left(u^{\alpha-p}-u_\l^{\alpha-p}\right) v^\beta \varphi_R^\theta(u^p-u_\l^p) \,dx \\ 
&+ \nu \alpha \int_{\hat B_{2R} \cap \{u\ge u_\l, \ v>v_\l\}} u_\l^{\alpha-p} \left(v^\beta-v_\lambda^\beta\right) \varphi_R^\theta (u^p-u_\l^p) \,dx =: \I_{5,1}+\I_{5,2},
\end{aligned}
\end{equation}
where in the last inequality we exploited the fact that our system is cooperative and that all the integrals are computed in support of $(u-u_\l)^+$.
We start by estimating $\I_{5,1}$, i.e. thanks to the Mean Value Theorem we get
$$\begin{aligned}
\I_{5,1}&\le \begin{cases}
0   & \text{if}\quad \alpha-p\le 0 \\
\nu \alpha p (\alpha -p) \int_{\hat B_{2R}} u_\l^{\alpha-p-1} u^{p-1}v^\beta {(\xi_\lambda^+)^2} \varphi_R^\theta \, dx   & \text{if}\quad 0<\alpha-p<1 \\
\nu \alpha p \int_{\hat B_{2R}} u^{p-1}v^\beta(\xi_\lambda^+)^2 \varphi_R^\theta \, dx    & \text{if}\quad\alpha-p=1 \\
\nu \alpha p(\alpha-p) \int_{\hat B_{2R}} u^{\alpha-2}v^\beta(\xi_\lambda^+)^2 \varphi_R^\theta  \,dx  & \text{if}\quad\alpha-p>1
\end{cases}
\end{aligned}$$
The case $0<\alpha-p<1$ needs of another manipulation. Using \eqref{eq:ratiobdd} and again $c^u_2 \leq |x|^{\mu_2}u \leq C^u_2$ and  $c^v_2 \leq |x|^{\mu_2}v \leq C^v_2$, we get
\begin{equation*}\label{eq:I51_1}
\begin{aligned}
    \I_{5,1}&= \nu \alpha p (\alpha -p) \int_{\hat B_{2R}} \left(\frac{u}{u_\l}\right)^{\alpha-p-1} u^{\alpha-p-1} u^{p-1}v^\beta {(\xi_\lambda^+)^2} \varphi_R^\theta \, dx\\
    & \le \nu \alpha p (\alpha -p) (C_2^v)^\beta \int_{\hat B_{2R}} u^{\alpha-2}  \frac{1}{|x|^{\mu_2\beta}} {(\xi_\lambda^+)^2} \varphi_R^\theta \, dx\\
    & \le \nu \alpha \check C \int_{\hat B_{2R}}  u^{p^*-2} (\xi_\lambda^+)^2 \varphi_R^\theta \,dx.
\end{aligned}
\end{equation*}
Thus, we have
\begin{equation}\label{eq:I51_3}
    \begin{aligned}
\I_{5,1}&\le \begin{cases}
0   & \text{if}\quad \alpha-p\le 0 \\
\nu \alpha \check{\bar C} \int_{\hat B_{2R}} u^{p^\star-2} (\xi_\lambda^+)^2 \varphi_R^\theta  \,dx  & \text{if}\quad\alpha-p>0.
\end{cases}
\end{aligned}
\end{equation}
From now on, for $p>2$, we repeat verbatim the computation from \eqref{I4} to \eqref{stimaI4pmag2} and we get  
\begin{equation}\label{I51finalpmag2}
\begin{aligned}
    \mathcal{I}_{5,1} \le & \nu \alpha \check{\bar C} \int_{\hat B_{2R}} u^{p^\star-2} (\xi_\lambda^+)^2 \varphi_R^\theta  \,dx  \\
    \le & \frac{\nu \alpha \bar{\bar C}}{|\l|^{\kappa^*}} \int_{\hat{B}_{2R}\cap\{u\ge u_\l\}}\varphi_R^\theta u^2\(|\n u|+|\n u_\l|\)^{p-2} |\nabla \log u - \nabla \log u_\l|^2 \, dx + \frac{\nu \alpha \bar{\bar C}}{|\l|^{\kappa^*}R^\kappa}.
\end{aligned}
\end{equation}
The term $\mathcal{I}_{5,2}$ is much more standard. Thanks to the cooperativity of the system and the Mean Value Theorem we get
\begin{equation}\label{I52_0}
\begin{aligned}
    \mathcal{I}_{5,2} \le & \nu \alpha \int_{\hat B_{2R} \cap \{u\ge u_\l, \ v>v_\l\}} u_\l^{\alpha-p} \left(v^\beta-v_\lambda^\beta\right) \varphi_R^\theta (u^p-u_\l^p) \,dx \\
    \le & \nu \alpha \beta p \int_{\hat B_{2R} \cap \{u\ge u_\l, \ v>v_\l\}} \left(\frac{u_\l}{u}\right)^{\alpha-p} u^{\alpha-1} v^{\beta -1} \zeta_\lambda^+ \xi_\lambda^+ \varphi_R^\theta  \,dx \\
    \le & \nu \alpha \beta p C \int_{\hat B_{2R}}  u^{\alpha-1} v^{\beta -1} \zeta_\lambda^+ \xi_\lambda^+ \varphi_R^\theta  \,dx,
    \end{aligned}
\end{equation}
where in the last inequality we also used \eqref{eq:ratiobdd} in the case $\alpha-p<0$. Applying first \eqref{eq:estatInf} to the right hand side of \eqref{I52_0} and then Young's inequality  we get
\begin{equation}\label{I52_1}
\begin{aligned}
    \mathcal{I}_{5,2} 
    \le &  \nu \alpha \beta p C \int_{\hat B_{2R}}  \frac{1}{|x|^{\mu_2(p^*-2)}} \zeta_\lambda^+ \xi_\lambda^+ \varphi_R^\theta  \,dx,\\
    \le & \frac{\nu \alpha \beta p C}{2} \int_{\hat B_{2R}}  \frac{1}{|x|^{\mu_2(p^*-2)}} (\xi_\lambda^+)^2 \varphi_R^\theta  \,dx + \frac{\nu \alpha \beta p C}{2} \int_{\hat B_{2R}}  \frac{1}{|x|^{\mu_2(p^*-2)}}  (\zeta_\lambda^+)^2 \varphi_R^\theta  \,dx.
    \end{aligned}
\end{equation}
Now in order to estimate the right hand side of \eqref{I52_1} we have to repeat verbatim the computation from \eqref{I4} to \eqref{stimaI4pmag2}, for $p>2$, and thus we obtain
\begin{equation}\label{I52_1finalpmag2}
\begin{aligned}
\mathcal{I}_{5,2}  \le & \frac{\nu \alpha \beta p C}{2} \int_{\hat B_{2R}}  \frac{1}{|x|^{\mu_2(p^*-2)}} (\xi_\lambda^+)^2 \varphi_R^\theta  \,dx + \frac{\nu \alpha \beta p C}{2} \int_{\hat B_{2R}}  \frac{1}{|x|^{\mu_2(p^*-2)}}  (\zeta_\lambda^+)^2 \varphi_R^\theta  \,dx\\
\le & \frac{\nu \alpha \beta p C}{|\l|^{\kappa^*}} \int_{\hat{B}_{2R}\cap\{u\ge u_\l\}}\varphi_R^\theta u^2\(|\n u|+|\n u_\l|\)^{p-2} |\nabla \log u - \nabla \log u_\l|^2 \, dx \\
&+ \frac{\nu \alpha \beta p C}{|\l|^{\kappa^*}} \int_{\hat{B}_{2R}\cap\{v\ge v_\l\}} \varphi_R^\theta v^2 \(|\n v|+|\n v_\l|\)^{p-2} |\nabla \log v - \nabla \log v_\l|^2 \, dx + \frac{2\nu \alpha \beta p C}{|\l|^{\kappa^*}R^\kappa}.
\end{aligned}
\end{equation}
Finally, collecting the estimates \eqref{I1}, \eqref{I2finale}, \eqref{I3finale},   \eqref{stimaI4pmag2}, \eqref{I5_0}, \eqref{I51finalpmag2}, \eqref{I52_1finalpmag2},   and using it in  \eqref{inizio}, we are able to deduce that
\begin{equation}\label{IequazioneSTimaFinale}
\begin{aligned}
&\(c_1-\frac{C+\nu \alpha (\bar{\bar C}+ \beta p C)}{|\l|^{\kappa^*}}\) \int_{\hat B_{2R}\cap \{u\ge u_\l\}} \varphi_R^\theta u^2\(|\n u|+|\n u_\l|\)^{p-2}|\n\log u -\n\log u_\l|^2 \, dx \\
&\le  \frac{\nu \alpha \beta p \hat C}{|\l|^{\kappa^*}} \int_{\hat{B}_{2R}\cap\{v\ge v_\l\}} \varphi_R^\theta v^2 \(|\n v|+|\n v_\l|\)^{p-2} |\nabla \log v - \nabla \log v_\l|^2 \, dx \\
& \quad + \frac{C}{R^\kappa} + \frac{C}{R^\frac{\kappa}{p}}  +\frac{C}{|\l|^{\kappa^*}R^\kappa} + \frac{\nu \alpha \bar{\bar C}}{|\lambda|^{\kappa^*}R^\kappa} +\frac{2\nu \alpha \beta p C}{|\l|^{\kappa^*}R^\kappa}.
\end{aligned}
\end{equation}
At this point we have to repeat the same arguments for the second equation \eqref{iniziosim2} and hence we have to give an estimate for all the terms $\E_1, \E_2, \E_3, \E_4, \E_5$. Thus we get
\begin{equation}\label{IIequazioneSTimaFinale}
\begin{aligned}
&\(c_2-\frac{C+\nu \beta (\tilde{C}+ \alpha p C)}{|\l|^{\kappa^*}}\) \int_{\hat B_{2R}\cap \{v\ge v_\l\}} \varphi_R^\theta v^2\(|\n v|+|\n v_\l|\)^{p-2}|\n\log v -\n\log v_\l|^2 \, dx \\
&\le  \frac{\nu \alpha \beta p \check C}{|\l|^{\kappa^*}} \int_{\hat{B}_{2R}\cap\{u\ge u_\l\}} \varphi_R^\theta v^2 \(|\n u|+|\n u_\l|\)^{p-2} |\nabla \log u - \nabla \log u_\l|^2 \, dx \\
& \quad + \frac{C}{R^\kappa} + \frac{C}{R^\frac{\kappa}{p}}  +\frac{C}{|\l|^{\kappa^*}R^\kappa} + \frac{\nu \beta \tilde{C}}{|\lambda|^{\kappa^*}R^\kappa} +\frac{2\nu \alpha \beta p C}{|\l|^{\kappa^*}R^\kappa}.
\end{aligned}
\end{equation}
Summing the two contribution \eqref{IequazioneSTimaFinale} and \eqref{IIequazioneSTimaFinale} we deduce
\begin{equation*}
    \begin{aligned}
        &\(c_1-\frac{C+\nu \alpha (\bar{\bar C}+ \beta p (C+\check C))}{|\l|^{\kappa^*}}\) \int_{\hat B_{2R}\cap \{u\ge u_\l\}} \varphi_R^\theta u^2\(|\n u|+|\n u_\l|\)^{p-2}|\n\log u -\n\log u_\l|^2 \, dx  \\
        & \ + \(c_2-\frac{C+\nu \beta (\tilde{C}+ \alpha p (C+\hat C))}{|\l|^{\kappa^*}}\) \int_{\hat B_{2R}\cap \{v\ge v_\l\}} \varphi_R^\theta v^2\(|\n v|+|\n v_\l|\)^{p-2}|\n\log v -\n\log v_\l|^2 \, dx \\
        & \le \frac{C}{R^\kappa} + \frac{C}{R^\frac{\kappa}{p}}  +\frac{C}{|\l|^{\kappa^*}R^\kappa} + \frac{\nu \beta \tilde{C}}{|\lambda|^{\kappa^*}R^\kappa} +\frac{2\nu \alpha \beta p C}{|\l|^{\kappa^*}R^\kappa}.
    \end{aligned}
\end{equation*}
For $|\l|$ sufficiently large, as $R$ goes to $+\infty$, we deduce that 
\begin{equation*}
\begin{split}
&\int_{\Sigma'_\l \cap \{u\ge u_\l\}}u^2\( |\nabla u|+|\n u_\l|\)^{p-2} |\n \log u- \n \log u_\l|^2  \, dx\\ 
& \ +  \int_{\Sigma'_\l \cap \{v\ge v_\l\}}  v^2\(|\n v|+|\n v_\l|\)^{p-2}|\n\log v -\n\log v_\l|^2 \, dx \le 0.
\end{split}
\end{equation*}
Now we have to estimate the right hand side of \eqref{stimaI4}, \eqref{eq:I51_3} and \eqref{I52_1} in the case $1<p<2$. We first remark that $2{ \vartheta}<0$ (for $n>2$) and, since $|x|\ge |x_\l|$, one has that $|x|^{2{\vartheta}}\le |x_\l|^{2{ \vartheta}}$.
Then
\begin{equation}\label{stimaI4pmin2}\begin{aligned}
\I_4 \le & \frac{C}{|\l|^{\kappa^*}}\int_{\hat{B}_{2R}\cap\{u\ge u_\l\}}|x_\l|^{2{ \vartheta}}\varphi_R^\theta |\nabla \log u -\nabla \log u_\l|^2 \, dx \\
&+\frac{C}{|\l|^{\kappa^*}} \int_{\hat{B}_{2R}\cap\{u\ge u_\l\}} |x|^{2{ \vartheta}} \(\log u - \log u_\l\)^2|\n \varphi_R|^2 \,dx.
\end{aligned}
\end{equation}
Let us consider $\check R=\max\{R_1, R_2\}$ and let $A_{\check R, \check R_0}=\overline{B_{\check R}(0_\l)\setminus B_{\check R_0}(0_\l)}$. Then we get 
$$\hat B_{2R}= \hat A_{\check R, \check R_0}\cup \(\hat B_{2R}\setminus \hat A_{\check R, \check R_0}\).$$
Thanks to Theorem \ref{thm:asymptEst} and in particular to \eqref{eq:estatInf}, we deduce that
\begin{equation}\label{stimaI4pmin2due}
\begin{aligned}
&\int_{\hat B_{2R}\setminus \hat A_{\check R, \check R_0}} |x_\l|^{2{ \vartheta}}\varphi_R^\theta |\nabla \log u - \n \log u_\l|^2 \, dx
\\
&{ =} C \int_{\hat B_{2R}\setminus \hat A_{\check R, \check R_0}} |x_\l|^{(\mu_2+1)(2-p)}|x_\l|^{-2\mu_2}\varphi_R^\theta |\nabla \log u - \n \log u_\l|^2 \, dx
\\
&\le C \int_{\hat B_{2R}\setminus \hat A_{\check R, \check R_0}} u_\l^2 \varphi_R^\theta\frac{ |\nabla \log u - \n \log u_\l|^2}{\(|\n u|+|\n u_\l|\)^{2-p}} \,dx.
 \end{aligned}
 \end{equation}
In the annulus $ A_{\check R, \check R_0}$ it is true that $|x_\l| \ge \check R_0$ (thus $|x_\l|^{2{ \vartheta}} \le \check R_0^{2{\vartheta}})$ and, since we are far from the reflected origin $0_\l$, we also obtain that $|\n u_\l|$ is bounded. Let $\mathcal L_u :=\displaystyle \inf_{B_{\check R}(0)\setminus B_{\check R_0}(0)} u$. Hence we get (by using \eqref{eq:ratiobdd} and the fact that $\(|\n u|+|\n u_\l|\)^{2-p}\le C$ far from $\{0, 0_\l\}$)
\begin{equation}\label{stimaI4pmin2tre}
\begin{aligned}
\int_{A_{\check R, \check R_0}} |x_\l|^{2{ \vartheta}}\varphi_R^\theta |\nabla \log u - \n \log u_\l|^2 \, dx & \le C\check R_0^{2\theta}\int_{A_{\check R, \check R_0}}\varphi_R^\theta |\nabla \log u - \n \log u_\l|^2 \, dx\\
&\le \frac{C\check R_0^{2\theta}}{\mathcal L_u^2}\int_{A_{\check R, \check R_0}} u_\l^2 \varphi_R^\theta \frac{ |\nabla \log u - \n \log u_\l|^2}{\(|\n u|+|\n u_\l|\)^{2-p}} \, dx. 
\end{aligned}
\end{equation}
Combining the two estimate \eqref{stimaI4pmin2due} and \eqref{stimaI4pmin2tre} in the first term of \eqref{stimaI4pmin2} and arguing as in \eqref{stimaI4pmag2} for the second term of \eqref{stimaI4pmin2} one yields to 
\begin{equation}\label{stimaI4pmin2unoprimo}
\begin{aligned}
&\I_4 \le \frac{C}{|\l|^{\kappa^*}} \int_{\hat{B}_{2R}\cap \{u\ge u_\l\}}u_\l^2 \varphi_R^\theta \frac{|\nabla \log u - \n\log u_\l|^2}{\(|\n u|+|\n u_\l|\)^{2-p}} \, dx + \frac{C}{|\l|^{\kappa^*}R^\kappa}. \end{aligned}
\end{equation}
Now, we go back to $\I_5$ and we want to estimate this term in the case $1<p<2$. Thanks to \eqref{eq:I51_3} and \eqref{I52_1}, we get
\begin{equation*}
\begin{aligned}
    \I_5 = \I_{5,1} +\I_{5,2}
    \le &   \nu \alpha \check{\bar C} \int_{\hat B_{2R}} u^{p^\star-2} (\xi_\lambda^+)^2 \varphi_R^\theta  \,dx +\frac{\nu \alpha \beta p C}{2} \int_{\hat B_{2R}}  \frac{1}{|x|^{\mu_2(p^*-2)}} (\xi_\lambda^+)^2 \varphi_R^\theta  \,dx \\
    &+ \frac{\nu \alpha \beta p C}{2} \int_{\hat B_{2R}}  \frac{1}{|x|^{\mu_2(p^*-2)}}  (\zeta_\lambda^+)^2 \varphi_R^\theta  \,dx\\
    \le & K_1 \int_{\hat B_{2R}}  \frac{1}{|x|^{\mu_2(p^*-2)}} (\xi_\lambda^+)^2 \varphi_R^\theta  \,dx + K_2\int_{\hat B_{2R}}  \frac{1}{|x|^{\mu_2(p^*-2)}}  (\zeta_\lambda^+)^2 \varphi_R^\theta  \, dx.
\end{aligned}
\end{equation*}
Repeating verbatim the argument used for $\I_4$, from \eqref{I4} to \eqref{stimaI4}, we get
\begin{equation}\label{I5pmin2_1}
    \begin{aligned}
        \I_5 \le &  \frac{C}{|\l|^{\kappa^*}} \int_{\hat B_{2R} \cap \{u\ge u_\l\}} |x|^{2\vartheta}|\nabla(\varphi_R^{\frac \theta 2}(\log u-\log u_\l)^+)|^2 \, dx\\
        &+ \frac{C}{|\l|^{\kappa^*}} \int_{\hat B_{2R} \cap \{v\ge v_\l\}} |x|^{2\vartheta}|\nabla(\varphi_R^{\frac \theta 2}(\log v-\log v_\l)^+)|^2 \, dx.
    \end{aligned}
    \end{equation}
Thus we can argue repeating verbatim the argument used for $\I_4$ in the case $1<p<2$, exactly from \eqref{stimaI4pmin2} to \eqref{stimaI4pmin2unoprimo}, for both the term in the right hand side of \eqref{I5pmin2_1}
\begin{equation} \label{finaleI5pmin2}
    \begin{aligned}
        \I_5 \le & \frac{\check C_1}{|\l|^{\kappa^*}} \int_{\hat{B}_{2R}\cap \{u\ge u_\l\}}u_\l^2 \varphi_R^\theta \frac{|\nabla \log u - \n\log u_\l|^2}{\(|\n u|+|\n u_\l|\)^{2-p}} \, dx \\
        &+ \frac{\check C_2}{|\l|^{\kappa^*}} \int_{\hat{B}_{2R}\cap \{v\ge v_\l\}}v_\l^2 \varphi_R^\theta \frac{|\nabla \log v - \n\log v_\l|^2}{\(|\n v|+|\n v_\l|\)^{2-p}} \, dx + \frac{C}{|\l|^{\kappa^*}R^\kappa}.
    \end{aligned}
\end{equation}
Hence, by collecting \eqref{I1bis}, \eqref{I2finale}, \eqref{I3finale},  \eqref{stimaI4pmin2unoprimo} and  \eqref{finaleI5pmin2} in \eqref{inizio}, we get that for $1<p<2$ it holds
\begin{equation}\label{IequazioneSTimaFinalepmin2}
    \begin{aligned}
    &\(C_p-\frac{C+\check C_1 K_1}{|\l|^{\kappa^*}}\) \int_{\hat B_{2R}\cap \{u\ge u_\l\}} u_\l^2\varphi_R^\theta \frac{ |\nabla \log u - \n \log u_\l|^2}{\(|\n u|+|\n u_\l|\)^{2-p}}\\
    &\le  \frac{\check C_2}{|\l|^{\kappa^*}} \int_{\hat{B}_{2R}\cap \{v\ge v_\l\}}v_\l^2 \varphi_R^\theta \frac{|\nabla \log v - \n\log v_\l|^2}{\(|\n v|+|\n v_\l|\)^{2-p}} \, dx + \frac {C}{ R^\kappa} + \frac{C}{R^\frac{\kappa}{p}}+ \frac{C}{|\l|^{\kappa^*}R^\kappa} .
    \end{aligned}
\end{equation}
As pointed out in the case $p>2$ we need to repeat the same arguments for the second equation \eqref{iniziosim2} and hence we have to give an estimate for all the terms $\E_1, \E_2, \E_3, \E_4, \E_5$. In particular, we are able to deduce that
\begin{equation}\label{IIequazioneSTimaFinalepmin2}
    \begin{aligned}
    &\(C_p-\frac{C+\check C_2 K_2}{|\l|^{\kappa^*}}\) \int_{\hat B_{2R}\cap \{v\ge v_\l\}} \varphi_R^\theta v_\l^2\frac{ |\nabla \log v - \n \log v_\l|^2}{\(|\n v|+|\n v_\l|\)^{2-p}} \, dx \\
    &\le  \frac{\check C_1}{|\l|^{\kappa^*}} \int_{\hat{B}_{2R}\cap \{u\ge u_\l\}}u_\l^2 \varphi_R^\theta \frac{|\nabla \log u - \n\log u_\l|^2}{\(|\n u|+|\n u_\l|\)^{2-p}} \, dx + \frac {C}{R^\kappa} + \frac{C}{R^\frac{\kappa}{p}}+ \frac{C}{|\l|^{\kappa^*}R^\kappa} .
    \end{aligned}
\end{equation}
Summing up the two contributions of \eqref{IequazioneSTimaFinalepmin2} and \eqref{IIequazioneSTimaFinalepmin2}, we obtain
\begin{equation*}\label{STimaFinalepmin2}
    \begin{aligned}
    \(C_p-\frac{\bar C_1}{|\l|^{\kappa^*}}\)& \int_{\hat B_{2R}\cap \{u\ge u_\l\}} \varphi_R^\theta u_\l^2\frac{ |\nabla \log u - \n \log u_\l|^2}{\(|\n u|+|\n u_\l|\)^{2-p}}\\
   & + \(C_p-\frac{\bar C_2}{|\l|^{\kappa^*}}\) \int_{\hat B_{2R}\cap \{v\ge v_\l\}} \varphi_R^\theta v_\l^2\frac{ |\nabla \log v - \n \log v_\l|^2}{\(|\n v|+|\n v_\l|\)^{2-p}}\\
    &\le  \frac{\check C_1}{|\l|^{\kappa^*}}  \frac {C}{R^\kappa} + \frac{C}{R^\frac{\kappa}{p}}+ \frac{C}{|\l|^{\kappa^*}R^\kappa}.
    \end{aligned}
\end{equation*}
Once again we can choose $|\l|$ large enough so that, as $R$ goes to $+\infty$, it yields
\begin{equation*}
\int_{\Sigma'_\l \cap \{u\ge u_\l\}}u_\l ^2\frac{ |\nabla \log u - \n \log u_\l|^2}{\(|\n u|+|\n u_\l|\)^{2-p}} \, dx + \int_{\Sigma'_\l \cap \{ v \ge v_\l\}} v_\l ^2\frac{ |\nabla \log v - \n \log v_\l|^2}{\(|\n v|+|\n v_\l|\)^{2-p}} \, dx \le 0
\end{equation*}
Hence, in both cases, $\log u -\log u_\l$ and $\log v -\log v_\l$ are constants and, since $\log u -\log u_\l =0$ and $\log v -\log v_\l =0$ on $T_\l$, then $\log u-\log u_\l =0$ and $\log v -\log v_\l =0$ in $\Sigma'_\l \cap \{u\ge u_\l\}$. Therefore we deduce both $u\le u_\l$ and $v\le v_\l$ in $\Sigma_\l$. Thus we proved that the set $\L_0$ is not empty and also that $\lambda_0$ exists and is finite.

\vspace{0.5cm}

{\bf Step 2:} 	$\lambda_0=0$. \\
The second part of the moving plane procedure is based on a contradiction argument, i.e.~assuming that $\lambda_0 < 0$.  Arguing as in the proof of \textbf{Step 1}  we will obtain a contradiction  proving that, simultaneously $u\le u_{\lambda_0 +\varepsilon}$ and $v\le v_{\lambda_0 +\varepsilon}$ in $\Sigma_{\lambda_0+\varepsilon}$ for all $0 < \varepsilon \le \overline{\varepsilon}$ for some $\overline{\varepsilon}>0$.

A key ingredient that we will use in the rest of the proof is the well known classical \emph{strong comparison principle}. In order to exploit this result, we notice that
from \textbf{Step 1} and by continuity it holds that
$$
u\le u_{\lambda_0} \quad \text{and} \quad v\le v_{\lambda_0}\quad \text{in}\quad\Sigma_{\lambda_0}\,.
$$
Let us first define the critical sets
$$
\mathcal{Z}_u:=\{x \in \R^n \ \text{s.t.} \ |\nabla u(x)|=0 \} \text{ and } \mathcal{Z}_v:=\{x \in \R^n \ \text{s.t.} \ |\nabla v(x)|=0 \}. 
$$

By \cite[Theorem 1.4]{damascelli} we deduce that
$u\equiv u_{\lambda_0}$ or $u<u_{\lambda_0}$  in any connected component $\mathcal C_u$ of $\Sigma_{\lambda_0}\setminus{\mathcal Z_u}$, and $v\equiv v_{\lambda_0}$ or $v<v_{\lambda_0}$ in any connected component $\mathcal C_v$ of $\Sigma_{\lambda_0}\setminus{\mathcal Z_v}$.
We will frequently use the fact that $\mathcal Z_u$ and $\mathcal Z_v$ have zero  Lebesgue measure \cite{DS1, EMM}.\\

\noindent We first assume that $\Sigma_{\lambda_0}\setminus{\mathcal Z_u}$ and $\Sigma_{\lambda_0}\setminus{\mathcal Z_v}$  have only one connected component.  We observe that the cases $u\equiv u_{\lambda_0}$ and $v\equiv v_{\lambda_0}$ are not possible since, by Theorem \ref{thm:asymptEst} and more precisely by \eqref{eq:estatInf}, there exists $B_{\check R_0}(0_{\lambda_0})$ where $u<u_{\lambda_0}$ and $v<v_{\lambda_0}$; this fact immediately implies that $u<u_{\lambda_0}$ in $\Sigma_{\lambda_0}\setminus \mathcal Z_u$ and $v<v_{\lambda_0}$ in $\Sigma_{\lambda_0}\setminus \mathcal Z_v$.\\

Let us assume that there exist at least two connected components of $\Sigma_{\lambda_0}\setminus{\mathcal Z_u}$ and two connected components of  $\Sigma_{\lambda_0}\setminus{\mathcal Z_v}$.  Theorem \ref{stimagradiente2} implies that $\mathcal Z_u$ and $\mathcal Z_v$ are bounded. As a consequence, we have that only one connected component of the set $\Sigma_{\lambda_0}\setminus{\mathcal Z_u}$ and one of the set $\Sigma_{\lambda_0}\setminus{\mathcal Z_v}$ can be unbounded. We refer to both these unbounded connected components as $\mathcal C_{1,u}$ and $\mathcal C_{1,v}$ and set
$$
\mathcal C_\lambda^u \,:=\,(\mathcal C_{1,u}^c\cap \Sigma_{\lambda_0})\cup R_\lambda(\mathcal C_{1,u}^c\cap \Sigma_{\lambda_0}) \text{ and } \mathcal C_\lambda^v \,:=\,(\mathcal C_{1,v}^c\cap \Sigma_{\lambda_0})\cup R_\lambda(\mathcal C_{1,v}^c\cap \Sigma_{\lambda_0}).
$$
If $u\equiv u_{\lambda_0}$ in $\mathcal C_{1,u}$ and $v\equiv v_{\lambda_0}$ in $\mathcal C_{1,v}$ it is easy to see that, by symmetry,  $\mathcal C^u_\lambda$ contains at least one connected component of $\mathbb{R}^n\setminus \mathcal Z_u$ and $\mathcal C^v_\lambda$ contains at least one connected component of $\mathbb{R}^n\setminus \mathcal Z_v$. But this is impossible as proved in \cite[Theorem 1.4]{DS1}, \cite[Lemma 5]{MMPS} and \cite[Proof of Theorem 1.1]{EMM}.
If else $u\equiv u_{\lambda_0}$ in $\mathcal C_{2,u}$ for some bounded component $\mathcal C_{2,u}$ and $v\equiv v_{\lambda_0}$ in $\mathcal C_{2,v}$ for some bounded component $\mathcal C_{2,v}$, then in this case we set
$$
\mathcal C_\lambda^u\,:=\,\mathcal C_{2,u} \cup R_\lambda(\mathcal C_{2,u}) \text{ and } \mathcal C_\lambda^v\,:=\,\mathcal C_{2,v} \cup R_\lambda(\mathcal C_{2,v})\,,
$$
 and also in this case, by symmetry,  $\mathcal C^u_\lambda$ would contain at least one connected component of $\mathbb{R}^n\setminus \mathcal Z_u$ and $\mathcal C^v_\lambda$  would contain at least one connected component of $\mathbb{R}^n\setminus \mathcal Z_v$ thus providing a contradiction. Resuming we just proved that
\begin{equation}\nonumber
u< u_{\lambda_0}\quad\text{in}\quad
\Sigma_{\lambda_0}\setminus{\mathcal Z_u} \quad \text{ and }  \quad v< v_{\lambda_0}\quad\text{in}\quad
\Sigma_{\lambda_0}\setminus{\mathcal Z_v}\,.
\end{equation}
 Now, recalling that $\mathcal Z_u$ and $\mathcal Z_v$ are bounded by Theorem \ref{stimagradiente2},  we fix $\overline{R} >0$ in such a way that
 $$
 \mathcal Z_u, \mathcal Z_v\subset B_{\overline{R}}(0)\,,
 $$
 and, for $\tau>0$, we let $\mathcal{Z}_u^\tau$ be an open set containing $\mathcal Z_u$ such that $\mathcal L (\mathcal Z_u^\tau)<\tau$ and $\mathcal{Z}_v^\tau$ be an open set containing $\mathcal Z_v$ such that $\mathcal L (\mathcal Z_v^\tau)<\tau$ (that exists since $\mathcal L (\mathcal Z_u)=\mathcal L (\mathcal Z_v)=0$).
 Then, for  $\delta,\varepsilon, \overline{R},\tau >0$, we denote by
 \begin{itemize}
     \item[(i)] $B_{\overline{R},\varepsilon}:= B^c_{\overline{R}}(0) \cap \Sigma_{\lambda_0+\varepsilon}$;

     \item[(ii)]  $S_{\delta,\varepsilon}:= \left((\Sigma_{\lambda_0+\varepsilon}\setminus\Sigma_{\lambda_0-\delta})\cap B_{\overline{R}}(0)\right)\cup (\mathcal Z_u^\tau\cap \Sigma_{\lambda_0-\delta}) \cup (\mathcal Z_v^\tau\cap \Sigma_{\lambda_0-\delta})$;

     \item[(iii)]   $K_\delta:= \overline{B_{\overline{R}}(0)\cap \Sigma_{\lambda_0-\delta}}\cap (\mathcal Z_u^\tau \cup \mathcal Z_v^\tau)^c$,
 \end{itemize}
where $\delta\le \overline{\delta}$ so that $K_\delta$ is nonempty.
We remark that this construction gives $$\Sigma_{\lambda_0 + \varepsilon} = B_{\overline{R},\varepsilon} \cup S_{\delta,\varepsilon} \cup K_\delta.$$
We also remark that, since $K_\delta$ is compact, then by the uniform continuity of $u$, $u_\l$ and $v$, $v_\l$, for $\overline{\varepsilon}>0$ small enough one has that $u<u_{{\lambda_0}+\varepsilon}$  and $v<v_{{\lambda_0}+\varepsilon}$ in $K_\delta$ for every $\varepsilon\le \overline{\varepsilon}$. Moreover, we point out the existence of $\tilde{R_0}$ such that $u<u_{\lambda_0 + \varepsilon}$ and  $v<v_{\lambda_0 + \varepsilon}$ in $B_{\tilde{R_0}}(0_{\lambda_0 + \varepsilon}) \subset \Sigma_{\lambda_0 + \varepsilon}$ for every $\varepsilon \le \overline{\varepsilon}$, and with $\tilde{R_0}$ independent of $\varepsilon$ as done in \textbf{Step 1}. \\

In the sequel, for $R>\overline{R}$, we consider $\varphi_R \in \mathcal C^\infty_c(B_{2R}(0))$ a cut-off function with $0\le \varphi_R \le 1$, $\varphi_R \equiv 1$ on $B_R(0)$ and $|\nabla \varphi_R|\le \frac{2}{R}$.

Then, letting $\theta>\max\{2, p\}$, we consider the following test functions
\begin{equation*}\label{test2}
\begin{split}
\Phi_{1, \lambda_0+\varepsilon}=\varphi_R^\theta u^{1-p}\xi_{p, \lambda_0+\varepsilon}^+ \chi_{\Sigma_{\lambda_0+\varepsilon}},\qquad \Phi_{2, \lambda_0+\varepsilon}=\varphi_R^\theta u_{\lambda_0+\varepsilon}^{1-p} \xi_{p, \lambda_0+\varepsilon}^+ \chi_{\Sigma_{\lambda_0+\varepsilon}},\\
\Psi_{1, \lambda_0+\varepsilon} =\varphi_R^\theta v^{1-p}\zeta_{p, \lambda_0+\varepsilon}^+ \chi_{\Sigma_{\lambda_0+\varepsilon}},\qquad \Psi_{2, \lambda_0+\varepsilon} =\varphi_R^\theta v_{\lambda_0+\varepsilon}^{1-p}\zeta_{p, \lambda_0+\varepsilon}^+ \chi_{\Sigma_{\lambda_0+\varepsilon}}.
\end{split}
\end{equation*}
Arguing similarly to \textbf{Step 1}, let us take $\Phi_{1, \lambda_0 +\varepsilon}$ and $\Psi_{1, \lambda_0 +\varepsilon}$ as test functions in both equation of \eqref{P}, and $\Phi_{2, \lambda_0 +\varepsilon}$, $\Psi_{2, \lambda_0 +\varepsilon}$ in \eqref{pbriflessa} respectively and, reasoning as in \textbf{Step 1}, one yields to
\begin{equation}\label{finesimeq1}
\begin{aligned}
&\int_{\hat B_{2R}} \varphi_R^\theta \left(|\nabla u|^{p-2}\nabla u\cdot\nabla(u^{1-p}\xi_{p,\lambda_0+\varepsilon}^+)-|\nabla u_{\lambda_0+\varepsilon}|^{p-2} \nabla u_{\lambda_0+\varepsilon} \cdot \nabla(u_{\lambda_0+\varepsilon}^{1-p} \xi_{p,\lambda_0+\varepsilon}^+) \right) \, dx\\
& \quad \le-\theta\int_{\hat B_{2R}} \varphi_R^{\theta-1} u^{1-p} \xi_{p, {\lambda_0+\varepsilon}}^+ |\nabla u|^{p-2}\nabla u\cdot \nabla \varphi_R \, dx \\
&\qquad + \theta \int_{\hat B_{2R}} \varphi_R^{\theta-1} u_\l^{1-p} \xi_{p, \lambda_0+\varepsilon}^+|\nabla u_{\lambda_0+\varepsilon}|^{p-2}\nabla u_{\lambda_0+\varepsilon} \cdot \nabla \varphi_R \,dx \\
&\qquad +\int_{\hat B_{2R}} \left(u^{p^*-p}-u_{\lambda_0+\varepsilon}^{p^*-p} \, dx \right) \varphi_R^\theta\xi_{p, \lambda_0+\varepsilon}^+ \, dx \\
& \qquad + \nu \alpha \int_{\hat B_{2R}} \left(u^{\alpha-p} v^\beta-u_{\lambda_0+\varepsilon}^{\alpha-p} v_{\lambda_0+\varepsilon}^\beta \right)\varphi_R^\theta\xi_{p,\lambda_0+\varepsilon}^+ \, dx,
\end{aligned}
\end{equation}
and to
\begin{equation}\label{finesimeq2}
\begin{aligned}
&\int_{\hat B_{2R}} \varphi_R^\theta \left(|\nabla v|^{p-2}\nabla v \cdot \nabla(v^{1-p}\zeta_{p, \lambda_0+\varepsilon}^+) -|\nabla v_{\lambda_0+\varepsilon}|^{p-2} \nabla v_{\lambda_0+\varepsilon} \cdot \nabla(v_{\lambda_0+\varepsilon}^{1-p} \zeta_{p, \lambda_0+\varepsilon}^+)\right) \, dx\\
& \quad \le - \theta\int_{\hat B_{2R}}\varphi_R^{\theta-1} v^{1-p} \zeta_{p,\lambda_0+\varepsilon}^+ |\nabla v|^{p-2} \nabla v \cdot \nabla \varphi_R \, dx \\
& \qquad + \theta \int_{\hat B_{2R}} \varphi_R^{\theta-1} v_{\lambda_0+\varepsilon}^{1-p} \zeta_{p,\lambda_0+\varepsilon}^+ |\nabla v_{\lambda_0+\varepsilon}|^{p-2} \nabla v_{\lambda_0+\varepsilon} \cdot \nabla \varphi_R \,dx \\
& \qquad +\int_{\hat B_{2R}}\left(v^{p^*-p}-v_{\lambda_0+\varepsilon}^{p^*-p}\right) \varphi_R^\theta\zeta_{p,\lambda_0+\varepsilon}^+ \,dx \\
& \qquad + \nu \beta \int_{\hat B_{2R}}\left(u^\alpha v^{\beta-p}-u_{\lambda_0+\varepsilon}^\alpha v_{\lambda_0+\varepsilon}^{\beta-p} \right) \varphi_R^\theta\zeta_{p,\lambda_0+\varepsilon}^+\,dx.
\end{aligned}
\end{equation}
We focus our attention only on \eqref{finesimeq1}.  Arguing as in \textbf{Step 1}, in both cases $p>2$ and $1<p<2$, we get
\begin{equation}\label{stimastep2eq1}
\begin{aligned}
&c_1 \int_{\hat B_{2R} \cap \{u\ge u_{\lambda_0 +\varepsilon}\}} \varphi_R^\theta u^2\(|\n u|+|\n u_{\lambda_0 +\varepsilon}|\)^{p-2}|\n\log u -\n\log u_{\lambda_0 +\varepsilon}|^2 \, dx
\\
& \quad \le  \int_{\hat B_{2R} \cap B_{\overline{R},\varepsilon}} (u^{p^*-p}-u_{\lambda_0 +\varepsilon}^{p^*-p}) \varphi_R^\theta \xi_{p, \lambda_0+\varepsilon}^+ \, dx + \int_{\hat B_{2R}\cap S_{\delta,\varepsilon}} (u^{p^*-p}-u_{\lambda_0 +\varepsilon}^{p^*-p}) \varphi_R^\theta \xi_{p, \lambda_0+\varepsilon}^+ \, dx \\
&\qquad + \nu \alpha \int_{\hat B_{2R} \cap B_{\overline{R},\varepsilon}} \left(u^{\alpha-p} v^\beta-u_{\lambda_0+\varepsilon}^{\alpha-p} v_{\lambda_0+\varepsilon}^\beta \right)\varphi_R^\theta\xi_{p,\lambda_0+\varepsilon}^+ \, dx\\
& \qquad + \nu \alpha \int_{\hat B_{2R} \cap S_{\delta,\varepsilon}} \left(u^{\alpha-p} v^\beta-u_{\lambda_0+\varepsilon}^{\alpha-p} v_{\lambda_0+\varepsilon}^\beta \right)\varphi_R^\theta\xi_{p,\lambda_0+\varepsilon}^+ \, dx + \frac{C}{R^{\frac \kappa p}}+ \frac{C}{R^\kappa}.
\end{aligned}
\end{equation}
We point out that we have used once again the fact that $\frac{u_{\lambda_0 +\varepsilon}}{u}\ge \tilde c$ for every $0\le \varepsilon\le \bar\varepsilon$ as to deduce \eqref{stimaulu}. \\
In order to estimate the first term on the right hand side of \eqref{stimastep2eq1}, we argue exactly in the same way of estimating $\I_4$ in \textbf{Step 1} (taking into account Remark \ref{remexterior}), where here $\overline{R}$ plays the role of $\l$ in \textbf{Step 1}. Hence we get
\begin{equation*}\label{stima1step2}
\begin{aligned}
&\int_{\hat B_{2R}\cap B_{\overline{R},\varepsilon}}(u^{p^*-p}-u_{\lambda_0 +\varepsilon}^{p^*-p})\varphi_R^\theta\xi_{p,\lambda_0+\varepsilon}^+ \, dx \le \frac{C}{R^\kappa { \bar R^{\kappa^*}}}
\\
&+ \frac{C}{\overline{R}^{\kappa^*}} \int_{\hat B_{2R} \cap B_{\overline{R},\varepsilon} \cap \{u\ge u_{\lambda_0 +\varepsilon}\}}\varphi_R^\theta u^2\(|\n u|+|\n u_{\lambda_0 + \varepsilon}|\)^{p-2} |\n\log u -\n\log u_{\lambda_0 +\varepsilon}|^2 \, dx.
\end{aligned}
\end{equation*}
For the second term on the right hand side of \eqref{stimastep2eq1} we argue as in \textbf{Step 1}, obtaining
\begin{equation}\label{stima2step21}
\begin{aligned}
&\int_{\hat B_{2R}\cap S_{\delta,\varepsilon}}(u^{p^*-p}-u_{\lambda_0 +\varepsilon}^{p^*-p})\varphi_R^\theta \xi_{p,\lambda_0+\varepsilon}^+ \,dx\le C_u\int_{\hat B_{2R}\cap S_{\delta,\varepsilon} \cap \{u\ge u_{\lambda_0 +\varepsilon}\}}(\log u -\log u_{\lambda_0 +\varepsilon})^2 \, dx,
\end{aligned}
\end{equation}
where $$C_u=\sup_{S_{\overline{\delta}, \overline{\varepsilon}}} u^{p^*-2}.$$
Now we need to split the estimate in two cases; if $p>2$ we apply a suitable weighted Poincar\'e inequality to the right hand side of \eqref{stima2step21} which can be found in \cite[Theorem 3.2]{DS1}. Hence in this case one has
\begin{equation}\label{stima2step2}
\begin{aligned}
&\int_{\hat B_{2R}\cap S_\delta^\varepsilon}(u^{p^*-p}-u_{\lambda_0 +\varepsilon}^{p^*-p})\varphi_R^\theta\xi_{p,\lambda}^+ \, dx
\\
&\le C^2_p(S_{\delta,\varepsilon})C_u\int_{\hat B_{2R}\cap S_{\delta,\varepsilon} \cap \{u\ge u_{\lambda_0 +\varepsilon}\}}|\n u|^{p-2}|\n\log u - \n \log u_{\lambda_0 +\varepsilon}|^2 \, dx
\\
&\le \frac{C^2_p(S_{\delta,\varepsilon} )C_u}{\displaystyle \inf_{S_{\overline{\delta}, \overline{\varepsilon}}} u^2} \int_{\hat B_{2R}\cap S_{\delta,\varepsilon} \cap \{u\ge u_{\lambda_0 +\varepsilon}\}}u^2\(|\n u|+|\n u_{\lambda_0 +\varepsilon}|\)^{p-2}|\n\log u - \n \log u_{\lambda_0 +\varepsilon}|^2 \, dx ,
\end{aligned}
\end{equation}
where $C_p(\Omega)$ is the Poincar\'e constant tends to zero as $|\Omega|\to 0$. If $1<p<2$ one can apply the classical Poincar\'e inequality in order to have
\begin{equation}\label{stima2bisstep2}
\begin{aligned}
&\int_{\hat B_{2R}\cap S_{\delta,\varepsilon}}(u^{p^*-p}-u_{\lambda_0 +\varepsilon}^{p^*-p})\varphi_R^\theta\xi_{p,\lambda_0+\varepsilon}^+ \, dx
\\
&\le C^2_p(S_{\delta,\varepsilon})C_u\int_{\hat B_{2R}\cap S_{\delta,\varepsilon} \cap \{u\ge u_{\lambda_0 +\varepsilon}\}}|\n\log u - \n \log u_{\lambda_0 +\varepsilon}|^2 \,dx
\\
&\le \frac{CC^2_p(S_{\delta,\varepsilon}) C_u}{\displaystyle \inf_{S_{\overline{\delta}, \overline{\varepsilon}}}u^2} \int_{\hat B_{2R}\cap S_{\delta,\varepsilon} \cap \{u\ge u_{\lambda_0 +\varepsilon}\}}u^2\(|\n u|+|\n u_{\lambda_0 +\varepsilon}|\)^{p-2}|\n\log u - \n \log u_{\lambda_0 +\varepsilon}|^2 \, dx,
\end{aligned}
\end{equation}
which can be deduced since in $\Sigma_{\lambda_0 + \varepsilon} \setminus B_{\tilde{R}_0}(0_{\lambda_0 + \varepsilon})$ one has that
$$\(|\n u|+|\n u_\l|\)^{2-p}\le C,$$
for some constant $C$ that does not depend on $\varepsilon \le \overline{\varepsilon}$.

Going back to \eqref{stimastep2eq1} we have to give an estimate for the third and the fourth term on the right hand side. Starting by the third term we can argue in the same way of $\I_5$
\begin{equation*}
    \begin{aligned}
    & \int_{\hat B_{2R} \cap B_{\overline{R},\varepsilon}} \left(u^{\alpha-p} v^\beta-u_{\lambda_0+\varepsilon}^{\alpha-p} v_{\lambda_0+\varepsilon}^\beta \right)\varphi_R^\theta\xi_{p,\lambda_0+\varepsilon}^+ \, dx\\
    \le &  \frac{C}{\overline{R}^{\kappa^*}} \int_{\hat{B}_{2R}\cap  B_{\overline{R},\varepsilon} \cap \{u\ge u_{\lambda_0+\varepsilon}\}} \varphi_R^\theta u^2\(|\n u|+|\n u_{\lambda_0+\varepsilon}|\)^{p-2} |\nabla \log u - \nabla \log u_{\lambda_0+\varepsilon}|^2 \, dx \\
&+ \frac{C}{\overline{R}^{\kappa^*}} \int_{\hat{B}_{2R}  \cap B_{\overline{R},\varepsilon} \cap\{v\ge v_{\lambda_0+\varepsilon}\}} \varphi_R^\theta v^2 \(|\n v|+|\n v_{\lambda_0+\varepsilon}|\)^{p-2} |\nabla \log v - \nabla \log v_{\lambda_0+\varepsilon}|^2 \, dx + \frac{C}{R^\kappa{ \bar R^{\kappa^*}}}.
\end{aligned}
\end{equation*}
For the last term of \eqref{stimastep2eq1}, we argue as in \textbf{Step 1}, so that
\begin{equation*}
\begin{aligned}
    &\int_{\hat B_{2R} \cap S_{\delta,\varepsilon}} \left(u^{\alpha-p} v^\beta-u_{\lambda_0+\varepsilon}^{\alpha-p} v_{\lambda_0+\varepsilon}^\beta \right)\varphi_R^\theta\xi_{p,\lambda_0+\varepsilon}^+ \, dx \\
    & \le C \int_{\hat B_{2R} \cap S_{\delta,\varepsilon}} u^{p^\star-2} (\xi_\lambda^+)^2 \varphi_R^\theta  \,dx  + C \int_{\hat B_{2R} \cap S_{\delta,\varepsilon}}  v^{p^*-2} (\zeta_\lambda^+)^2 \varphi_R^\theta  \,dx\\
    &\le C_u\int_{\hat B_{2R}\cap S_{\delta,\varepsilon} \cap \{u\ge u_{\lambda_0 +\varepsilon}\}}(\log u -\log u_{\lambda_0 +\varepsilon})^2 \, dx +C_v\int_{\hat B_{2R}\cap S_{\delta,\varepsilon} \cap \{v\ge v_{\lambda_0 +\varepsilon}\}}(\log v -\log v_{\lambda_0 +\varepsilon})^2 \, dx,
\end{aligned}
\end{equation*}
where
$$C_u:= \sup_{S_{\overline{\delta}, \overline{\varepsilon}}} u^{p^*-2} \qquad \text{and} \qquad  
C_v:=\sup_{S_{\overline{\delta}, \overline{\varepsilon}}} v^{p^*-2}.$$
At this point we can argue in two different cases $p>2$ as in \eqref{stima2step2} and $1<p<2$ as in \eqref{stima2bisstep2}.
\begin{equation*}
    \begin{aligned}
        &\int_{\hat B_{2R} \cap S_{\delta,\varepsilon}} \left(u^{\alpha-p} v^\beta-u_{\lambda_0+\varepsilon}^{\alpha-p} v_{\lambda_0+\varepsilon}^\beta \right)\varphi_R^\theta\xi_{p,\lambda_0+\varepsilon}^+ \, dx \\
        \le& \frac{CC^2_p(S_{\delta,\varepsilon}) C_u}{\displaystyle \inf_{S_{\overline{\delta}, \overline{\varepsilon}}}u^2} \int_{\hat B_{2R}\cap S_{\delta,\varepsilon} \cap \{u\ge u_{\lambda_0 +\varepsilon}\}}u^2\(|\n u|+|\n u_{\lambda_0 +\varepsilon}|\)^{p-2}|\n\log u - \n \log u_{\lambda_0 +\varepsilon}|^2 \, dx \\
        &+ \frac{CC^2_p(S_{\delta,\varepsilon}) C_v}{\displaystyle \inf_{S_{\overline{\delta}, \overline{\varepsilon}}}v^2} \int_{\hat B_{2R}\cap S_{\delta,\varepsilon} \cap \{v\ge v_{\lambda_0 +\varepsilon}\}}v^2\(|\n v|+|\n v_{\lambda_0 +\varepsilon}|\)^{p-2}|\n\log v - \n \log v_{\lambda_0 +\varepsilon}|^2 \, dx.
    \end{aligned}
\end{equation*}
Analogously, we can obtain the same estimate starting from \eqref{finesimeq2} and repeating verbatim the above argument. At the end, we have to collect all the estimates { in \eqref{stimastep2eq1}} and sum the contributions, in order to deduce that
\begin{equation*}\label{stima3step2}
\begin{aligned}
&\check c \int_{\hat B_{2R} \cap \{u\ge u_{\lambda_0 +\varepsilon}\}} \varphi_R^\theta u^2\(|\n u|+|\n u_{\lambda_0 +\varepsilon}|\)^{p-2}|\n\log u -\n\log u_{\lambda_0 +\varepsilon}|^2 \, dx
\\
& \quad  + \check c \int_{\hat B_{2R} \cap \{v\ge v_{\lambda_0 +\varepsilon}\}} \varphi_R^\theta v^2\(|\n v|+|\n v_{\lambda_0 +\varepsilon}|\)^{p-2}|\n\log v -\n\log v_{\lambda_0 +\varepsilon}|^2 \, dx
\\
&\le   \frac{C}{R^{\frac \kappa p}}+\frac{C}{R^\kappa}+\frac{C}{\bar R^{\kappa^*}}\displaystyle \int_{(\hat B_{2R}\cap B_{\bar R, \varepsilon})\cap \{u\ge u_{\lambda_0 +\varepsilon}\}}\varphi_R^\theta u^2\(|\n u|+|\n u_{\lambda_0 +\varepsilon}|\)^{p-2}|\n\log u -\n\log u_{\lambda_0 +\varepsilon}|^2 \, dx \\ 
&\quad +\frac{\tilde C C^2_p(S_\delta^\varepsilon)C_u}{\displaystyle \inf_{S_{\overline{\delta}}^{\overline{\varepsilon}}}u^2}\int_{S_\delta^\varepsilon\cap \{u\ge u_{\lambda_0 +\varepsilon}\}}{u^2}\(|\n u|+|\n u_{\lambda_0 +\varepsilon}|\)^{p-2}|\n\log u - \n \log u_{\lambda_0 +\varepsilon}|^2 \, dx\\
& \quad +\frac{C}{\bar R^{\kappa^*}} \int_{(\hat B_{2R}\cap B_{\bar R, \varepsilon})\cap \{v\ge v_{\lambda_0 +\varepsilon}\}}\varphi_R^\theta v^2\(|\n v|+|\n v_{\lambda_0 +\varepsilon}|\)^{p-2}|\n\log v -\n\log v_{\lambda_0 +\varepsilon}|^2 \, dx \\ 
& \quad +\frac{\bar C C^2_p(S_\delta^\varepsilon)C_v}{\displaystyle \inf_{S_{\overline{\delta}}^{\overline{\varepsilon}}}v^2}\int_{S_\delta^\varepsilon\cap \{v\ge v_{\lambda_0 +\varepsilon}\}}{v^2}\(|\n v|+|\n v_{\lambda_0 +\varepsilon}|\)^{p-2}|\n\log v - \n \log v_{\lambda_0 +\varepsilon}|^2 \, dx,
\end{aligned}
\end{equation*}
where $\check c= \min\{c_1,c_2\}$.
Now we take care of the variable parameters $\bar R, \delta, \bar\varepsilon$. First we fix $\bar R$ large such that  
$$\frac{C}{\check c \bar R^{\kappa^*}}<1.$$
Then, since $C_p^2(\Omega)$ goes to zero if the Lebesgue measure of $\Omega$ goes to zero, we choose $\delta,\bar\varepsilon,\tau$ small so that 
$$
\frac{CC^2_p(S_\delta^\varepsilon)C_u}{\check c \displaystyle \inf_{S_{\overline{\delta}}^{\overline{\varepsilon}}}u^2}<1 \quad \text{ and } \quad \frac{CC^2_p(S_\delta^\varepsilon)C_v}{\check c \displaystyle \inf_{S_{\overline{\delta}}^{\overline{\varepsilon}}}v^2}<1
$$ 
for every $0\le \varepsilon \le \bar\varepsilon$.
Hence it follows that
\begin{equation*}\label{stima4step2}
\begin{aligned}
&\int_{\hat B_{2R}\cap \{u\ge u_{\lambda_0 +\varepsilon}\}}u^2\(|\n u|+|\n u_{\lambda_0 +\varepsilon}|\)^{p-2}|\n\log u -\n\log u_{\lambda_0 +\varepsilon}|^2 \, dx\\
&+ \int_{\hat B_{2R}\cap \{v\ge v_{\lambda_0 +\varepsilon}\}}v^2\(|\n v|+|\n v_{\lambda_0 +\varepsilon}|\)^{p-2}|\n\log v -\n\log v_{\lambda_0 +\varepsilon}|^2 \, dx \le \frac{C}{R^{\frac \kappa p}}+\frac{C}{R^\kappa},
\end{aligned}
\end{equation*}
getting again (as $R\to+\infty$)
\begin{equation*}\label{stima6step2}
\begin{aligned}
&\int_{\Sigma_{\lambda_0 + \varepsilon} \cap \{u\ge u_{\lambda_0 +\varepsilon}\}}u^2\(|\n u|+|\n u_{\lambda_0 +\varepsilon}|\)^{p-2}|\n\log u -\n\log u_{\lambda_0 +\varepsilon}|^2 \, dx\\
&+ \int_{\hat B_{2R}\cap \{v\ge v_{\lambda_0 +\varepsilon}\}}v^2\(|\n v|+|\n v_{\lambda_0 +\varepsilon}|\)^{p-2}|\n\log v -\n\log v_{\lambda_0 +\varepsilon}|^2 \, dx= 0,
\end{aligned}
\end{equation*}
which gives that $u\le u_{\lambda_0 +\varepsilon}$ and $v\le v_{\lambda_0 +\varepsilon}$ in $\Sigma_{\lambda_0 + \varepsilon}$. This contradicts the definition of $\lambda_0$ and it proves that $\lambda_0=0$. 

\vspace{0.5cm}

\textbf{Step 3:} Conclusion.\\
The symmetry of the solution $(u,v)$ follows by performing the moving plane procedure in the opposite direction, and hence this gives the symmetry of $u$ along the $e_1$-direction. Repeating the same arguments with respect to each direction $\eta\in \mathbb S^{n-1}_+$, then one deduces that $(u,v)$ has radial symmetry {with respect to the} origin and that is a radially decreasing function. \end{proof}

\end{document}